\magnification 1200
\input amstex.tex
\input amsppt.sty

\vsize24.8truecm
\hsize15.2truecm
\voffset-5truemm
\hoffset+3truemm


\author{Maxim Nazarov}\endauthor
\title{Yangian of the Queer Lie Superalgebra}\endtitle


\expandafter\ifx\csname maxim.def\endcsname\relax \else\endinput\fi
\expandafter\edef\csname maxim.def\endcsname{%
 \catcode`\noexpand\@=\the\catcode`\@\space}
\catcode`\@=11

\mathsurround 1.6pt
\font\Bbf=cmbx12 

\def\hcor#1{\advance\hoffset by #1}
\def\vcor#1{\advance\voffset by #1}
\let\bls\baselineskip  \let\ignore\ignorespaces
\def\vsk#1>{\vskip#1\bls} \let\adv\advance 
\def\vv#1>{\vadjust{\vsk#1>}\ignore} \def\vvv#1>{\vadjust{\vskip#1}\ignore}
\def\vvn#1>{\vadjust{\nobreak\vsk#1>\nobreak}\ignore}
\def\vvvn#1>{\vadjust{\nobreak\vskip#1\nobreak}\ignore}
\def\emph#1{{\it #1\/}}

 \let\nt\noindent 
\def\nn#1>{\noalign{\vskip #1pt}} \def\NN#1>{\openup#1pt}
 
\let\Sum\sum \def\sum{\Sum\limits} 
\let\Prod\prod \def\prod{\Prod\limits} \let\Int\int \def\int{\Int\limits}

\let\=\m@th \def\&{.\kern.1em} \def\>{\!\;} \def\:{\!\!\;}

\ifx\plainfootnote\undefined \let\plainfootnote\footnote \fi
\expandafter\ifx\csname amsppt.sty\endcsname\relax
 
\else \fi

\newbox\s@ctb@x
\def\s@ct#1 #2\par{\removelastskip\vsk>
 \vtop{\bf\setbox\s@ctb@x\hbox{#1} \parindent\wd\s@ctb@x
 \ifdim\parindent>0pt\adv\parindent.5em\fi\item{#1}#2\strut}%
 \nointerlineskip\nobreak\vtop{\strut}\nobreak\vsk-.4>\nobreak}

\newbox\t@stb@x
\def\gadv{\global\advance} \def\gad#1{\gadv#1 1} 
\def\l@b@l#1#2{\def\n@@{\csname #2no\endcsname}%
 \if *#1\gad\n@@ \expandafter\xdef\csname @#1@#2@\endcsname{\the\Sno.\the\n@@}%
 \else\expandafter\ifx\csname @#1@#2@\endcsname\relax\gad\n@@
 \expandafter\xdef\csname @#1@#2@\endcsname{\the\Sno.\the\n@@}\fi\fi}
\def\l@bel#1#2{\l@b@l{#1}{#2}\?#1@#2?}
\def\?#1?{\csname @#1@\endcsname}
\def\[#1]{\def\n@xt@{\ifx\t@st *\def\n@xt####1{{\setbox\t@stb@x\hbox{\?#1@F?}%
 \ifnum\wd\t@stb@x=0 {\bf???}\else\?#1@F?\fi}}\else
 \def\n@xt{{\setbox\t@stb@x\hbox{\?#1@L?}\ifnum\wd\t@stb@x=0 {\bf???}\else
 \?#1@L?\fi}}\fi\n@xt}\futurelet\t@st\n@xt@}
\def\(#1){{\rm\setbox\t@stb@x\hbox{\?#1@F?}\ifnum\wd\t@stb@x=0 ({\bf???})\else
 (\?#1@F?)\fi}}
\def\dff{\expandafter\d@f} \def\d@f{\expandafter\def}
\def\edff{\expandafter\ed@f} \def\ed@f{\expandafter\edef}

\newcount\Sno \newcount\Lno \newcount\Fno
\def\Section#1{\gad\Sno\Fno=0\Lno=0\s@ct{\the\Sno.} {#1}\par} \let\Sect\Section
\def\section#1{\gad\Sno\Fno=0\Lno=0\s@ct{} {#1}\par} \let\sect\section
\def\l@F#1{\l@bel{#1}F} \def\<#1>{\l@b@l{#1}F} \def\l@L#1{\l@bel{#1}L}
\def\Tag#1{\tag\l@F{#1}} \def\Tagg#1{\tag"\llap{\rm(\l@F{#1})}"}
\def\Th#1{Theorem \l@L{#1}} \def\Lm#1{Lemma \l@L{#1}}
\def\Prop#1{Proposition \l@L{#1}}
\def\Cr#1{Corollary \l@L{#1}} \def\Cj#1{Conjecture \l@L{#1}}
 
\def\Proof#1.{\demo{\it Proof #1}}

 \def\setparindent{\edef\Parindent{\the\parindent}}
\def\Appendix{\Sno=64\let\p@r@\z@ 
\def\Section##1{\gad\Sno\Fno=0\Lno=0 \s@ct{} \hskip\p@r@ Appendix \char\the\Sno
 \if *##1\relax\else {.\enspace##1}\fi\par} \let\Sect\Section
\def\section##1{\gad\Sno\Fno=0\Lno=0 \s@ct{} \hskip\p@r@ Appendix%
 \if *##1\relax\else {.\enspace##1}\fi\par} \let\sect\section
\def\l@b@l##1##2{\def\n@@{\csname ##2no\endcsname}%
 \if *##1\gad\n@@
\expandafter\xdef\csname @##1@##2@\endcsname{\char\the\Sno.\the\n@@}%
\else\expandafter\ifx\csname @##1@##2@\endcsname\relax\gad\n@@
 \expandafter\xdef\csname @##1@##2@\endcsname{\char\the\Sno.\the\n@@}\fi\fi}}

\let\logo@\relax
\let\m@k@h@@d\makeheadline \let\m@k@f@@t\makefootline
\def\makeheadline{\ifnum\pageno=1\headline={\hfil}\fi\m@k@h@@d}
\def\makefootline{\ifnum\pageno=1\footline={\hfil}\fi\m@k@f@@t}


\def\1{^{-1}}

\def\a{{\frak{a}}}
\def\An{{\Cal A_{\ts n}}}
\def\Anp{{\Cal A_{\ts n^\prime}}}
\def\Anpn{{\Cal A_{\ts n+ n^\prime}}}

\def\al{\sigma}

\def\b{\frak{b}}

\def\bi{{\bar\imath}}
\def\bj{{\bar\jmath}}
\def\bk{{\bar k}}
\def\bl{{\bar l}}
\def\br{{\bar r}}

\def\CC{{\Bbb C}}
\def\CN{\CC^{N|N}}

\def\de{\delta}
\def\De{\Delta}
\def\deg{{\operatorname{deg}\ts}}
\def\Den{{\Delta^{\hskip-4pt\o}}}
\def\DYN{\Cal D\!\YN}

\def\End{\operatorname{End}\hskip1pt}
\def\EndCN{\End(\CN)}
\def\enddemos{{\ $\square$\enddemo}}
\def\ep{\varepsilon}
\def\et{\eta\ts}

\def\g{\frak{g}}
\def\ga{\gamma}
\def\ge{\geqslant}
\def\glN{\frak{gl}_{N|N}}
\def\gr{\operatorname{gr}\!}
\def\grm{\gr\hskip.5pt\raise-3pt\hbox{-}\hskip-2.5pt}

\def\h{\frak{h}}

\def\id{{\operatorname{id}}}
\def\io{\iota}

\def\lc{{\ts,\hskip.95pt\ldots\ts,\ts}}
\def\le{\leqslant}

\def\M{{\Cal M}}

\def\o{{\raise.24ex\hbox{$\sssize\kern.1em\circ\kern.1em$}}}

\def\om{\omega}
\def\ot{\otimes}

\def\ph{\varphi}

\def\qN{\frak{q}_N}

\def\R{\Cal R}

\def\S{{\Cal S}}

\def\Sh{\operatorname{S}(\h)}

\def\TB{T^{\hskip.5pt\ast\hskip-.5pt}}
\def\ts{{\hskip1pt}}

\def\Ua{\operatorname{U}(\a)}
\def\Ug{\operatorname{U}(\g)}
\def\UglNu{\operatorname{U}(\hskip.5pt\glN[u]\ts)}
\def\Ugp{\operatorname{U}(\g^\prime)}
\def\Uh{\operatorname{U}(\h)}
\def\UN{\operatorname{U}(\qN)}

\def\up{\alpha}

\def\Yh{\operatorname{Y}(\qN,h)}
\def\YN{\operatorname{Y}(\qN)}

\def\YNP{\operatorname{Y}^{\ast\hskip-1pt}(\qN)}
\def\YNPs{\gr_s\!\YNP}
\def\YNs{\gr_s\YN}
\def\YNS{\operatorname{Y}^{\hskip.5pt\prime\hskip-.5pt}(\qN)}

\def\ZZ{\Bbb Z}


\line{\Bbf Yangian of the Queer Lie Superalgebra\hfill}
\bigskip\bigskip
\line{\bf Maxim Nazarov\hfill}
\bigskip
\line{Department of Mathematics, University of York,
York YO1 5DD, England\ts.\hfill}
\line{E-mail: \tt mln1{\@}york.ac.uk\hfill}
\bigskip\smallskip


\noindent
{\bf Abstract.} Consider the complex matrix Lie superalgebra $\glN$ with the
standard generators $E_{ij}$ where $i,j=\pm1\lc\!\!\pm\!N\!$.
Define an involutory automorphism $\et$ of $\glN$ by
$\et(E_{ij})=E_{-i,-j}$. The twisted polynomial current Lie superalgebra
$$
\g=\{\,X(u)\in\glN[u]\,:\,\et(X(u))=X(-u)\,\}
$$
has a natural Lie co-superalgebra structure. We quantise
the universal enveloping algebra $\Ug$ as a 
co\ts-Poisson Hopf superalgebra. For the quantised algebra
we give a description of the centre,
and construct the double in the sense of Drinfeld.
We also construct a wide class of irreducible representations of
the quantised algebra.


\smallskip
\section{1. Introduction}

\noindent 
In this article we will work with certain Lie superalgebras [K]
over the complex field $\CC$. 
Their universal enveloping algebras are
$\ZZ_2$-graded associative unital algebras, and we
will always keep to the following convention.
Let $\operatorname{A}$ and $\operatorname{B}$ be any two associative
complex $\ZZ_2$-graded algebras. Their tensor product
$\operatorname{A}\ot\operatorname{B}$ will be
a $\ZZ_2$-graded algebra such that for any homogeneous
$X,X^\prime\in\operatorname{A}$ and $Y,Y^\prime\in\operatorname{B}$
$$
\align
(X\ot Y)\ts (X^\prime\ot Y^\prime)&=X\ts X^\prime\ot Y\ts Y^\prime
\cdot(-1)^{\ts\deg X^\prime\deg Y},
\\
\deg\ts(X\ot Y)&=\deg X+\deg Y\ts.
\endalign
$$
Throughout this article we will denote by $\theta$ the isomorphism
$\operatorname{A}\ot\operatorname{B}\to\operatorname{B}\ot\operatorname{A}$
defined by
$$
X\ot Y\mapsto Y\ot X\cdot(-1)^{\ts\deg X\deg  Y}.
$$
If the algebra $\operatorname{A}$ is unital denote by $\iota_p$
its embedding into the tensor product $\operatorname{A}^{\!\ot n}$
as the $p$-\ts th tensor factor:
$$
\iota_p(X)=1^{\ot\ts (p-1)}\ot X\ot1^{\ot\ts(n-p)}\,,
\qquad
1\le p\le n.
$$
We will also use various embeddings of the algebra
$\operatorname{A}^{\!\ot\ts m}$ into $\operatorname{A}^{\!\ot\ts n}$
for any $m\leqslant n$.
For any choice of pairwise distinct indices $p_1\lc p_m\in\{\ts1\lc n\ts\}$
and an element $X\in\operatorname{A}^{\!\ot m}$ of the form
$X=X^{(1)}\ot\ldots\ot X^{(m)}$ we will denote
$$
X_{p_1\ldots p_m}=\ts
\iota_{p_1}\bigl(X^{(1)}\bigr)\ts\ldots\,\ts\iota_{p_m}\bigl(X^{(m)}\bigr)
\in\operatorname{A}^{\!\ot n}.
$$

Let $\a$ be an arbitrary finite-dimensional Lie superalgebra.
Then consider the {\it polynomial current\/} Lie superalgebra $\a[u]$. 
It consists of the polynomial functions of a complex variable $u$
valued in $\a$. For any two such functions their supercommutator in
$\a[u]$ is determined pointwise.
Let $K\in\a^{\otimes2}$ be an $\a$-invariant element: we have the equality
$
\bigl[\ts X_1+X_2\ts,K\ts\bigr]=0
$
in $\a^{\otimes2}$ for any $X\in\a\ts$. Here $\a$ is regarded as a subspace
in the enveloping algebra
$\Ua$ and the square brackets stand for
the supercommutator. Also suppose that
$K$ is of $\ZZ_2$-degree zero and symmetric\ts: $K_{12}=K_{21}$.
Then the rational function
$$
r(u,v)=\frac K{u-v}
\Tag{1.1}
$$
of two complex variables $u\,,v$ satisfies  
the {\it classical Yang\ts-Baxter equation} for $\a[u]\ts$:
$$
\bigl[\ts r_{12}(u,v),r_{13}(u,w)\ts\bigr]+
\bigl[\ts r_{12}(u,v),r_{23}(v,w)\ts\bigr]+
\bigl[\ts r_{13}(u,w),r_{23}(v,w)\ts\bigr]=0
\Tag{1.2}
$$
in $\a^{\ot3}$.
This can be verified by momentary calculation. Furthermore, \text{the function}
\(1.1) is antisymmetric:
$$
r_{12}(u,v)+r_{21}(v,u)=0\,.
$$
Therefore the co-supercommutator  
$\ph:\ts\a[u]\to\a[u]^{\ot2}=\a^{\ot2}[u,v]$
\text{can be defined by}
$$
\ph\bigl(X(u)\bigr)=\bigl[\ts X_1(u)+X_2(v)\,,\ts r(u,v)\ts\bigr].
\Tag{1.3}
$$
This definition 
makes $\a[u]$ into a Lie bi-superalgebra.
In particular, if $\a$ is a simple Lie algebra and $K$ is the Casimir element,
one gets a natural Lie bialgebra structure on $\a[u]$.
It gives rise to a natural co-Poisson structure on the universal
enveloping algebra $\operatorname{U}(\a[u]\ts)$,
which is a co\ts-commutative Hopf algebra by definition. 
The more general case of a simple Lie superalgebra $\a$ was considered in [LS].

Now consider the queer Lie superalgebra $\qN$.
This is the most interesting super-analogue of the general linear Lie
algebra $\frak{gl}_N$, see for instance [S2].
We will realise $\qN$
as a subalgebra in the general linear Lie superalgebra $\glN$.
Let the indices $i,j$ run through $\pm\,1\lc\pm\,N$.
We will always write $\bi=0$ if $i>0$ and $\bi=1$ if $i<0$. Consider the
$\ZZ_2$-graded vector space $\CN$.
Let $e_i\in\CN$ be an element of the standard basis.
The $\ZZ_2$-gradation on $\CN$ is defined so that $\deg e_i=\bi$.
Let $E_{ij}\in\EndCN$ be the standard matrix units. The algebra $\EndCN$
is $\ZZ_2$-graded so that $\deg E_{ij}=\bi+\bj$. We will
also regard $E_{ij}$ as generators of the complex Lie superalgebra $\glN$.
The {\it queer\/} classical Lie superalgebra $\qN$ is
the fixed point subalgebra in $\glN$ with respect to
the involutive automorphism
$$
\et:\thinspace E_{ij}\mapsto E_{-i,-j}
\Tag{1.35}
$$

The queerness of $\qN$ reveals itself in that all the symmetric
$\qN$-invariants in $\qN^{\ot2}$ of $\ZZ_2$-degree zero are trivial: 
for $\a=\qN$ we always have $K\in\ts\CC\!\ts\cdot\!E^{\ts\ot2}$ where
$$
E=E_{11}+E_{-1,-1}+\ldots+E_{NN}+E_{-N,-N}\,.
\Tag{2.3.0}
$$
Hence in this case the co-supercommutator \(1.3) vanishes and there is no
natural Lie
bi-superalgebra structure on $\qN[u]$. However such a structure
can be defined, in compensation, on the
{\it twisted polynomial current\/} Lie superalgebra
$$
\g=\bigl\{\ts X(u)\in\glN[u]\ts:\,\et(X(u))=X(-u)\ts\bigr\}\,.
\Tag{1.0}
$$

Our definition is based on the following general scheme [A1,A2,FR].
Let $\a\ts,K$ be arbitrary as above and $\om$ be an automorphism of the Lie
superalgebra $\a$ of finite order $n$.
Let $\zeta$ be a primitive $n$-th root of unity. Generalising (1.1) put
$$
r(u,v)\ =\sum_{m\,\in\,\ZZ_n}\,\frac{\id\ot\om^m\,(K)}{u-\zeta^m\,v}\,.
\Tag{1.4}
$$
\proclaim{Proposition 1.1}
Suppose that $\om^{\ts\ot2}\,(K)=\zeta\thinspace K$.
Then the function \(1.4) is antisymmetric and obeys
the classical Yang\ts-Baxter equation \(1.2).
\endproclaim

\demo{Proof}
The function $r(u,v)$ determined by \(1.1) satisfies the equation \(1.2).
Let us apply to the left-hand side of \(1.2) with that $r(u,v)$
the operator $\id\ts\ot\ts\om^{\ts k}\ts\ot\ts\om^{\ts l}$ in $\a^{\ot3}$
and substitute $\zeta^k\,v\ts,\zeta^l\,w$ for $v\ts,w$ respectively.
Taking then the sum over $k,l\in\ZZ_n$ and using 
$\om^{\ts\ot2}\,(K)=\zeta\thinspace K$,
we will obtain
the left-hand side of \(1.2) with
the function $r(u,v)$ determined by \(1.4).
\text{For the latter function $r(u,v)$ we also have}
$$
\align
r_{21}(v,u)
&=
\sum_{m\,\in\,\ZZ_n}\,\frac{\om^m\ot\id\,(K_{21})}{v-\zeta^m\,u}\ =
\sum_{m\,\in\,\ZZ_n}\,\frac{\om^m\ot\id\,(K)}{v-\zeta^m\,u}
\\
&=
\sum_{m\,\in\,\ZZ_n}\,\frac{\id\ot\om^{-m}\,(K)}{\zeta^{-m}\,v-u}\ \,=
\,-\,r(u,v)
\quad\square
\endalign
$$
\enddemo

\noindent
Note that for a simple Lie algebra $\a$ always
$\om^{\ts\ot2}\,(K)=K$, and in compliance with [BD]
this construction does not give any new solutions of \(1.2)\ts.
Let $\a$ be the Lie superalgebra $\glN$ and $\om$ be the
involutive automorphism \(1.35). The element
$$
P\,=\,\sum_{ij}\thinspace E_{ij}\ot E_{ji}\cdot(-1)^{\,\bj}
\Tag{1.5}
$$
of $\glN^{\,\ot2}$ is symmetric and $\glN$-invariant.
Moreover, we have $\et^{\ts\ot2}\ts(P)=-\ts P$.
Due to Proposition 1.1 by setting
$K=P$ in \(1.4) we get an antisymmetric solution of the
Yang\ts-Baxter equation \(1.2). Therefore \(1.3) defines a co-supercommutator
$\ph:\ts\g\to{\g}^{\ot2}$.
Thus we obtain a Lie bi-superalgebra structure on $\g$.

For any simple finite-dimensional Lie algebra $\a$, quantisation of
the co\ts-Poisson Hopf algebra $\!\operatorname{U}(\a[u]\ts)\!$ was
described in [D1].
The quantised Hopf algebra is denoted by $\operatorname{Y}(\a)$
and called the Yangian of the Lie algebra $\a$. The algebra
$\operatorname{Y}(\a)$ contains the universal enveloping 
algebra $\Ua$ as a subalgebra.
However, the case $\a=\frak{sl}_N$ is exceptional since only for
$\a=\frak{sl}_N$ there exists a homomorphism
$\operatorname{Y}(\a)\to\Ua$ identical on the subalgebra $\Ua$, see
[\ts D1\ts,\ts Theorem 9\ts].
There is also a Hopf algebra $\operatorname{Y}(\frak{gl}_N)$,
which is a quantisation of
the co\ts-Poisson Hopf algebra $\operatorname{U}(\frak{gl}_N[u]\ts)$.
Again, the algebra
$\operatorname{Y}(\frak{gl}_N)$ contains the enveloping algebra
$\operatorname{U}(\frak{gl}_N)$ as a subalgebra, and admits
a homomorphism
$\operatorname{Y}(\frak{gl}_N)\to\operatorname{U}(\frak{gl}_N)$
identical on $\operatorname{U}(\frak{gl}_N)$. Moreover, the algebra
$\!\operatorname{Y}(\frak{gl}_N)\!$ can be defined entirely
in terms of the classical representation theory [O1].
For further details on the Yangian $\operatorname{Y}(\frak{gl}_N)$
see [MNO] and references therein.

\vbox{
The main aim of this article is to define the Yangian of
the Lie superalgebra $\qN$. It cannot be defined as a quantisation
of the enveloping algbra $\operatorname{U}(\qN[u]\ts)$, because the latter
Hopf superalgebra has no natural co\ts-Poisson structure. Instead of $\qN[u]$
we will consider the twisted polynomial current Lie superalgebra $\g$.
In Section 2 we define a certain Hopf superalgebra $\Yh$
over the field $\CC\ts[[h]]$ of the formal power series in $h$.
The quotient $\Yh\,/\,h\Yh$ is isomorphic to $\Ug$ as
a co\ts-Poisson Hopf superalgebra.
All specialisations of $\Yh$ at $h\in\CC\setminus\{0\}$ are 
isomorphic to each other as Hopf superalgebras. The specialisation at $h=1$
will be denoted by $\YN$ and called the Yangian of Lie superalgebra $\qN$.
Similarly to the Yangian $\operatorname{Y}(\frak{gl}_N)$,
the algebra $\YN$ contains the enveloping algebra $\UN$ as a subalgebra,
and admits a homomorphism $\YN\to\UN$ identical on $\UN$.
In Section~3 we describe the centre of the $\ZZ_2$-graded algebra $\YN$.
In Section~4 we construct the double of this Yangian in the sense of [D3].
In Section~5 we study an analogue for $\YN$ of the Drinfeld functor [D2]
for the Yangian $\operatorname{Y}(\frak{gl}_N)$.
}


\section{\bf2.\ Definition of the Yangian}

\noindent
In this section we introduce the {\it Yangian\/} of the Lie superalgebra
$\qN$. This is a complex associative
unital $\ZZ_2$-graded algebra $\YN$ with the countable
set of generators $T^{(s)}_{ij}$ where $s=1,2,\ts\ldots$ and
$i,j=\pm1,\ts\dots\ts,\pm N$.
The $\ZZ_2$-gradation on the algebra $\YN$
is determined by setting $\deg T_{ij}^{(s)}=\bi+\bj$ for $s\ge1$.
To write down defining relations for these
generators we will employ the formal series
$$
T_{ij}(u)=
\de_{ij}\cdot1+T_{ij}^{(1)}\ts u^{-1}+T_{ij}^{(2)}\ts u^{-2}+\ldots
\Tag{3.0}
$$
from $\YN[[u\1]]$. Then for all possible indices $i,j$
and $k,l$ we have the relations
$$
\align
(u^2-v^2)
\cdot
\bigl[\,T_{ij}(u&)\ts,T_{kl}(v)\ts\bigr]\cdot
(-1)^{\ts\bi\ts\bk\ts+\ts\bi\ts\bl\ts+\ts\bk\ts\bl}=
\Tag{3.1}
\\
(u+v)
\cdot
\bigl(\ts T_{kj}(u)\,T_{il}(v)&-T_{kj}(v)\,T_{il}(u)\bigr)\,\,\ts\ -
\\
\quad\qquad
(u-v)
\cdot
\bigl(\ts T_{-k,\ts j}(u)\,T_{-i,\ts l}(v)
&-T_{k,\ts -j}(v)\,T_{i,\ts -l}(u)\bigr)
\cdot
(-1)^{\ts\bk\ts+\ts\bl}
\endalign
$$
in $\YN((u\1,v\1))$.
The square brackets here stand for the supercommutator.
Moreover, for all possible indices $i,j$  we impose the relations
$$
T_{ij}(-\ts u)=T_{-i,\ts -j}(u).
\Tag{3.2}
$$

We will also use the following matrix form of the relations \(3.1).
Regard $E_{ij}$ as elements of the algebra $\EndCN$.
Combine all the series \(3.0) into the single element
$$
T(u)\,=\,\sum_{ij}\ts E_{ij}\ot T_{ij}(u)
$$
of the algebra $\EndCN\ot\YN\ts[[u^{-1}]]$.
For any positive integer $n$ and each $s=1\lc n$ we denote
$$
T_s(u)=\iota_s\ot\id\ts\bigl(T(u)\bigr)\,\in\,\EndCN^{\ot n}\ot\YN\,[[u\1]]\ts.
\Tag{3.22}
$$

Regard \(1.5) as an element of the algebra $\EndCN^{\ot2}$.
Consider the element
$$
J\,=\,\sum_{i}\ts E_{i,-i}\cdot{(-1)}^{\ts\bi}
\Tag{1.55555}
$$
of the algebra $\EndCN$, it has $\ZZ_2$-degree one.  
Note that the supercommutant of this element in $\EndCN$
coincides with the image of the defining representation $\qN\to\EndCN$.
Introduce the rational function
of two complex variables $u\,,v$
$$
\gather
R\ts(u\,,v)=1-\frac{P}{u-v}+\frac{P J_1 J_2}{u+v}\ts=
\Tag{3.333}
\\
1-\sum_{ij}\ts E_{ij}\ot E_{ji}\cdot\frac{{(-1)}^{\ts\bj}}{u-v}\ts
-\sum_{ij}\ts E_{ij}\ot E_{-j,\ts -i}\cdot\frac{{(-1)}^{\ts\bj}}{u+v}
\endgather
$$
valued in the algebra $\EndCN^{\ot2}$.
Then the relations \(3.1) can be rewritten as
$$
\bigl(R\ts(u,v)\ot1\bigr)\cdot\ts T_1(u)\ts T_2(v)=
T_2(v)\ts T_1(u)\cdot\bigl(R\ts(u,v)\ot1\bigr).
\Tag{3.3}
$$
Namely, after multiplying each side of \(3.3) by $u^2-v^2$ it becomes
a relation in the algebra
$$
\EndCN^{\ot2}\ot\YN\,((u\1,v\1))
$$
equivalent to the collection of all relations \(3.1).
Also note that the function \(3.333) satisfies the
{\it quantum Yang\ts-Baxter equation\/}
for the algebra $\EndCN^{\ot\ts 3}(u,v,w)$
$$
R_{12}(u,v)\ts\ts R_{13}(u,w)\ts\ts R_{23}(v,w)=
R_{23}(v,w)\ts\ts R_{13}(u,w)\ts\ts R_{12}(u,v)\,.
\Tag{3.6}
$$

Furthermore, consider \(1.35) as an automorphism of the algebra $\EndCN$.
The collection of all relations \(3.2) is equivalent to the single equation
$$
\et\ot\id\,\bigl(\ts T(u)\bigr)=T(-u)\ts.
\Tag{3.4}
$$
Observe that by the definition \(3.333) of $R\ts(u,v)$ we also have
in $\EndCN^{\ot\ts2}(u,v)$
$$
\align
\et\ot\id\,\bigl(\ts R\ts(u,v)\bigr)&=R\ts(-\ts u,v)\ts,
\Tag{3.45}
\\
\id\ot\et\,\bigl(\ts R\ts(u,v)\bigr)&=R\ts(u,-\ts v).
\Tag{3.46}
\endalign
$$
We call the function \(3.333) the {\it rational $R$-matrix\,}
for the Lie superalgebra $\qN$.

For any $i,j$ put $F_{ij}=E_{ij}+E_{-i,-j}$.
Then we have the equality $\et(F_{ij})=F_{ij}$ in $\EndCN$.
We will also regard $F_{ij}$ as generators of the universal
enveloping algebra $\UN$. Due to \(3.1) there is a homomorphism
$$
\ts\YN\to\UN:\,T_{ij}(u)\mapsto \de_{ij}-F_{ji}\,u^{-1}\cdot{(-1)}^{\ts\bj}\ts.
\Tag{3.5}
$$
The relations \(3.1),\(3.2) now imply that
the assignment
$$
F_{ji}\mapsto-\,T_{ij}^{(1)}\cdot{(-1)}^{\ts\bj}
\Tag{3.51}
$$
defines embedding of $\ZZ_2$-graded associative unital algebras $\UN\to\YN$.
The fact that the homomorphism \(3.51) has no kernel, is due to
Theorem 2.3 below.

The homomorphism \(3.5) is identical on the subalgebra $\UN$. It will be called
the {\it evaluation homomorphism} for the algebra $\YN$ and denoted by $\pi_N$.
The element $T(u)\!$ of 
the algebra $\EndCN\ot\YN[[u\1]]$ is invertible, we~put
$$
T(u)\1=\sum_{i,j}\ts E_{ij}\ot\widetilde T_{ij}(u).
$$
Then the relations \(3.3),\(3.4) along with the identity
$$
R(u,v)\ts R(-u,-v)=1-\frac1{(u-v)^2}-\frac1{(u+v)^2}
\Tag{3.52}
$$
imply that the assignment $T_{ij}(u)\mapsto\widetilde T_{ij}(-\ts u)$
determines an automorphism of the algebra $\YN$.
This automorphism is evidently involutive.

We will use two different ascending $\ZZ$-filtrations on the algebra $\YN$.
They are obtained by assigning to the generator $T_{ij}^{(s)}$
the degree $s$ or $s-1$ respectively. The corresponding
$\ZZ$-graded algebras will be denoted by $\gr\ts\YN$ and $\grm\YN$.
Let $G_{ij}^{\ts(s)}\in\grm\YN$ be the element
corresponding to the generator $T_{ij}^{(s)}\in\YN$. 
The algebra $\grm\YN$ inherits $\ZZ_2$-gradation from $\YN$
such that $\deg\ts G_{ij}^{\ts(s)}=\bi\ts+\ts\bj$.

Take the enveloping algebra $\Ug$
of the twisted current Lie superalgebra \(1.0).
The algebra $\Ug$ also has a natural $\ZZ_2$-gradation: the $\ZZ_2$-degree
of the element
$$
F_{ij}^{(s)}=E_{ij}\,u^{\ts s}+E_{-i,-j}\ts(-\ts u)^s  
\Tag{3.56}
$$
equals $\bi\ts+\bj$ for any $s\ge0$.
We have the following easy observation.

\proclaim{Proposition 2.1}
The assignment for every $s\ge0$
$$
F_{ji}^{(s)}\mapsto-\ts\,G_{ij}^{\ts(s+1)}\cdot{(-1)}^{\ts\bj}
\Tag{3.5566}
$$
determines a surjective homomorphism $\Ug\to\grm\hskip.5pt\YN$
of $\ZZ_2$-graded algebras.
\endproclaim

\demo{Proof}
The elements \(3.56) generate the algebra $\Ug$.
The defining relations for these generators can be written as
$$
\gather
\bigl[\ts F_{ji}^{\ts(s)},F_{lk}^{\ts(r)}\ts\bigr]=
\de_{il}\,F_{jk}^{\ts(s+r)}-\,
\de_{kj}\,F_{li}^{\ts(s+r)}\cdot
{(-1)}^{\ts(\ts\bi\ts+\ts\bj\ts)(\ts\bl\ts+\ts\bk\ts)}\,\ts+
\Tag{2.2.1}
\\
\qquad\qquad
\de_{i,-l}\,F_{-j,k}^{\ts(s+r)}\cdot{(-1)}^s\,-
\de_{-k,j}\,F_{l,-i}^{\ts(s+r)}\cdot
{(-1)}^{\ts(\ts\bi\ts+\ts\bj\ts)(\ts\bl\ts+\ts\bk\ts)\ts+\ts s}
\endgather
$$
for all $r,s\ge0$ and
$$
F_{-j,-i}^{\ts(s)}=(-1)^s\cdot F_{ji}^{\ts(s)}.
\Tag{2.2.2}
$$
On the other hand, by \(3.1) we obtain the relations in
the algebra $\grm\YN$
$$
\gather
(-1)^{\ts\bi\ts\bk\ts+\ts\bi\ts\bl\ts+\ts\bk\ts\bl}\cdot
\bigl[\ts G_{ij}^{\ts(s)},G_{kl}^{\ts(r)}\ts\bigr]\ts=\,
\de_{kj}\,G_{il}^{\ts(s+r-1)}-\ts
\de_{il}\,G_{kj}^{\ts(s+r-1)}\,\ts+
\\
\qquad\qquad
\de_{k,-j}\,G_{-i,\ts l}^{\ts(s+r-1)}\cdot(-1)^{\bk\ts+\ts\bl\ts+\ts s}
\,-
\de_{-i,l}\,G_{k,-j}^{\ts(s+r-1)}
\cdot(-1)^{\bk\ts+\ts\bl\ts+\ts s}
\endgather
$$
for $r,s\ge1$. Due to \(3.2)
$$
G_{-i,-j}^{\ts(s)}=(-1)^s\cdot G_{ij}^{\ts(s)}.
$$
Comparison of these relations to \(2.2.1) and \(2.2.2) shows
that \(3.5566) determines a homomorphism $\Ug\to\grm\hskip.5pt\YN$.
This homomorphism is surjective and preserves $\ZZ_2$-gradation
by definition
\enddemos

\noindent
There is a natural Hopf superalgebra structure on $\YN$. Due to \(3.3),\(3.4)
the comultiplication $\De:\YN\to\YN\ot\YN$ can be defined by
$$
T_{ij}(u)\,\mapsto\,\sum_{k}\,\,
T_{ik}(u)\ot T_{kj}(u)\cdot
{(-1)}^{\ts(\ts\bi\ts+\ts\bk\ts)(\ts\bj\ts+\ts\bk\ts)}
\Tag{3.7}
$$
where the tensor product is taken over the subalgebra $\CC[[u^{-1}]]$
in $\YN\ts[[u\1]]$ and the index $k$ runs through $\pm1\lc\pm N$.
The counit $\ep:\YN\to\CC$ is defined so that
$\ep\!:\ts T^{(-s)}_{ij}\mapsto0$ for every $s\ge1$.
Then the assignment $T_{ij}(u)\mapsto\widetilde T_{ij}(u)$
determines the antipodal map $S\!:\YN\to\YN$. It is an antiautomorphism
of the $\ZZ_2$-graded algebra $\YN$. Note that
$\YN$ contains $\UN$ as a Hopf sub-superalgebra:
by the definitions \(3.51) and \(3.7) for any $F\in\qN$ we have
$$
\De\ts(F)=F\ot1+1\ot F\ts,
\ \quad
\ep(F)=0\ts,
\ \quad
S(F)=-F\ts.
$$

The comultiplication \(3.7) on the $\ZZ_2$-graded algebra $\YN$
allows us to define for any $n=1,2,\ldots$ a representation
$\YN\to\EndCN^{\ot n}$ depending on $n$ arbitrary complex parameters
$z_1\lc z_n$. Indeed, by comparing 
\(3.6),\(3.45) to \(3.3),\(3.4)
respectively we obtain that for any $z\in\CC$ the assignment
$$
\EndCN\ot\YN\ts[[u^{-1}]]
\rightarrow
\EndCN^{\ot\ts2}\ts[[u^{-1}]]
\,:\ 
T(u)\mapsto R(u,z)
\Tag{3.66}
$$
determines a representation $\YN\to\EndCN$. More explicitly, we have
$$
T_{ij}^{(s+1)}\,\mapsto\,
-\ts\bigl(\ts E_{ji}\,z^{\ts s}+E_{-j,\ts -i}\ts(-z)^{\ts s}\ts\bigr)
\cdot{(-1)}^{\ts\bj}\,,
\qquad
s\ge0\ts.
\Tag{2.19}
$$
When $z=0$ this representation $\YN$ 
can be also obtained from the standard representation $\UN\to\EndCN$
by virtue of the evaluation homomorphism \(3.5).                         
Now for any $z_1\lc z_n\in\CC$ take the tensor product
of the representations \(3.66) of the algebra $\YN$ with $z=z_1\lc z_n$.
Due to \(3.7) the respective homomorphism $\YN\to\EndCN^{\ot n}$
is determined by the assignment
$$
\gather
\EndCN\ot\YN\ts[[u^{-1}]]
\ \rightarrow\ 
\EndCN^{\ot(n+1)}\ts[[u^{-1}]]\ts:\ \ 
\Tag{3.666666}
\\
T(u)
\,\mapsto\,
R_{12}(u,z_1)\ts\ldots\ts R_{1,n+1}(u,z_n).
\endgather
$$

\proclaim{Proposition 2.2}
Let the complex parameters $z_1\lc z_n$ and positive integer $n$ vary.
Then the kernels of all representations \(3.666666) of $\YN$
have \text{zero intersection.\!}
\endproclaim

\demo{Proof}
Take any finite linear combination of the products
$T_{i_1j_1}^{(s_1)}\!\ldots T_{i_mj_m}^{(s_m)}\in\YN$
with certain complex coefficients
$A_{\ts i_1j_1\ldots i_mj_m}^{\ts(s_1\ldots s_m)}$
where the indices $s_1\lc s_m\ge1$ and the number $m\ge0$ may vary.
Consider the image of this combination under the representation 
$\YN\to\EndCN^{\ot n}$ determined by \(3.666666);
it depends on $z_1\lc z_n$ polynomially. Take the terms of
this polynomial which have the maximal total degree in $z_1\lc z_n$.
Let $A$ the sum of these terms and $d$ be their degree.

Consider the ascending $\ZZ$-filtration on algebra $\YN$ where
the generator $T_{ij}^{(s)}$ with $s\ge1$ has degree $s-1$.
Equip the tensor product $\YN^{\ot n}$ with the ascending $\ZZ$-filtration
where the degree is the sum of the degrees on the tensor factors. Then by
the definition \(3.7) under the comultiplication $\YN\to\YN^{\ot n}$ we have
$$
T_{ij}^{(s)}\mapsto
\sum_{1\le r\le n}1^{\ot(r-1)}\ot T_{ij}^{(s)}\ot1^{\ot(n-r)}\ +\
\text{lower degree terms\,,}\quad
s\ge1\,.
$$
Therefore $A\in\EndCN^{\ot n}$ coincides with the image of the sum
$$
\sum_{s_1+\ldots+s_m=\ts d+m}
A_{\ts i_1j_1\ldots i_mj_m}^{\ts(s_1\ldots s_m)}\
F_{j_1i_1}^{\ts(s_1-1)}\ldots\,F_{j_mi_m}^{\ts(s_m-1)}
\cdot(-1)^{\,m\ts+\ts\bj_1\ts+\ts\ldots\ts+\ts\bj_m}\in\Ug
$$
under the tensor product of the evaluation representations
$$
\Ug\to\EndCN:F_{ij}^{(s)}\,\mapsto\ts 
E_{ij}\,z^{\ts s}+E_{-i,-j}\ts(-\ts z)^s\,,
\quad s\ge0
$$
at the points $z=z_1\lc z_n\in\CC$; see the definition \(3.56) of the
element $F_{ij}^{(s)}\in\g$,
and the formula \(2.19)
for the representation $\YN\to\EndCN$ corresponding to $z\in\CC$.
Due to Proposition 2.1 it now suffices to show that when $z_1\lc z_n\in\CC$
and the positive integer $n$ vary, the kernels
of the tensor products
of the evaluation representations
of the algebra $\Ug$ at $z=z_1\lc z_n\in\CC$
have zero intersection.
This will also imply that the homomorphism \(3.5566) is injective.

The algebra $\Ug$ is a subalgebra in the universal enveloping algebra
of the Lie superalgebra $\glN[u]$. 
We will show that the intersection of the kernels of all finite tensor
products of 
evaluation representations $\UglNu\to\EndCN$, is zero.
Denote by $\varpi_n$ the supersymmetrisation map in
the tensor product $(\hskip.5pt\glN[u]\ts)^{\ot\ts n}$ normalised so that
$\varpi_n^{\ts2}=\varpi_n$. 
We will identify the vector space
$(\hskip.5pt\glN[u]\ts)^{\ot\ts n}$ with $\glN^{\,\ot\ts n}[u_1\lc u_n]$ where
$u_1\lc u_n$ are independent complex variables.

\newpage 

The vector space $\glN$ is identified with $\EndCN$.
Choose any~linear~basis $X_1\lc X_{4N^2}$ in $\glN$ such that $X_1=E$
as in \(2.3.0).
The element $X_1\in\glN$
is then identified with the operator $1\in\EndCN$. Take
any finite non-zero linear combination of the elements
$$
(X_{a_1}u^{s_1})\ldots(X_{a_m}u^{s_m})\in\UglNu
\Tag{2.3.1}
$$
where the indices $s_1\lc s_m\ge0$ and the number $m\ge0$ may vary.
We assume that for every fixed $m$ the elements
$$
\varpi_m\bigl(X_{a_1}u^{s_1}\ot\ldots\ot X_{a_m}u^{s_m}\bigr)\,\in\,\ts
(\hskip.5pt\glN[u]\ts)^{\ot\ts m}
\,=\,\EndCN^{\ot m}\ts[\ts u_1\lc u_m]
$$
are linearly independent. Further, we will suppose that in every product
\(2.3.1) the indices
$a_1\lc a_p>1$ for certain $p\le m$ while $a_{p+1}=\ldots=a_m=1$.
We will also suppose that $s_{p+1}\lc s_q>0$ for some $q\ge p$,
while $s_{q+1}=\ldots=s_m=0$.

For any $n\ge p$ consider the tensor product $\nu$ of the evaluation 
representations of the algebra $\UglNu$ at $u_1\lc u_n\in\CC$.
Let us denote by $\operatorname{P}$
the subspace~in $\EndCN^{\ot n}$ spanned by the vectors
$X_{b_1}\ot\ldots\ot X_{b_n}$ where either the number of
indices $b_r>1$ is less than $p$, or $b_r=1$ for at least one $r\le p$.
The image of \(2.3.1) under $\nu$ is a polynomial in
$u_1\lc u_n$ valued in $\EndCN^{\ot n}$, of the form
$$
p\ts!\,\varpi_p\bigl(X_{a_1}u^{s_1}\ot\ldots\ot X_{a_p}u^{s_p}\bigr)
\ot1^{\ts\ot(n-p)}
\,\cdot\!\!\!\prod_{p<r\le q}\!\!
(u_1^{s_r}+\ldots+u_n^{s_r})
\cdot n^{\ts m-q}
\Tag{2.3.2}
$$
plus the terms valued in the subspace $\operatorname{P}\subset\EndCN^{\ot n}$.
Here the tensor factor
$\varpi_p\bigl(X_{a_1}u^{s_1}\ot\ldots\ot X_{a_p}u^{s_p}\bigr)$
is regarded as an element of $\EndCN^{\ot p}\ts[\ts u_1\lc u_p\ts]$
by identifying this algebra with $(\hskip.5pt\glN[u]\ts)^{\ot\ts p}$.

The numbers $p$ for various products \(2.3.1) from our linear combination
may differ. Take those products \(2.3.1) where the number $p$ is maximal.
For any $n\ge p$ the images of the remaining products under $\nu$ 
are polynomials in $u_1\lc u_n$ taking values
in the subspace $\operatorname{P}\subset\EndCN^{\ot n}$.
But a non-zero linear combination of the polynomials
\(2.3.2) with the maximal $p$, cannot vanish identically for all $n\ge p$
by Poincar\'e\,-Birkhoff\,-Witt theorem
[\ts MM,\ts Theorem 5.15\ts] for Lie superalgebras
\enddemos

\noindent
In the course of the proof of Proposition 2.2 we established that
the homomorphism \(3.5566) is injective.
Together with Proposition 2.1, this yields the following result. 

\proclaim{Theorem 2.3}
$\ZZ_2$-graded algebras $\Ug$ and $\grm\hskip.5pt\YN$
\text{are isomorphic via \(3.5566).}
\endproclaim

\noindent
Let us now return to the first $\ZZ$-filtration on the algebra $\YN$.
Let $t_{ij}^{\ts(s)}$ be the element of the algebra $\gr\ts\YN\!$ 
corresponding to the generator $T_{ij}^{\ts(s)}\!\in\!\YN$.
The algebra $\gr\ts\YN$ inherits $\ZZ_2$-gradation from $\YN$
such that $\deg\ts t_{ij}^{\ts(s)}=\bi\ts+\ts\bj$.

\proclaim{Corollary 2.4}
The algebra $\gr\ts\YN$ is supercommutative with
free generators $t_{ij}^{\ts(s)}$ and $t_{i,-j}^{\ts(s)}$
where $s=1\ts,2\ts,\ts\ldots$ and $i,j=1\lc N$.
\endproclaim

\demo{Proof}
The $\ZZ$-graded algebra $\gr\ts\YN$ is supercommutative due to the relations
\(3.1). Moreover, by \(3.2) for any $s\ge1$ we have the relation
$t_{-i,-j}^{\ts(s)}={(-1)}^s\, t_{ij}^{\ts(s)}$.
The supercommuting generators $t_{ij}^{\ts(s)}$ with
$i>0$ are free due to Theorem 2.3
\enddemos

To finish this section let us show that the Hopf superalgebra $\YN$ provides
a quantisation of the co\ts-Poisson Hopf superalgebra $\Ug$ in the sense of
[D1]. Let $h$ be a formal parameter. Take the tensor product
$\CC[[h]]\ot\YN$ where $h$ has $\ZZ_2$-degree zero.
Denote by $\Yh$ the unital
subalgebra in this tensor product, generated by
all the elements $H^{(s)}_{ij}=T^{(s)}_{ij}h^{s-1}$ with $s\ge1$.
Due to Theorem~2.3 an \text{isomorphism} of $\ZZ_2$-graded algebras
$\!\Yh/\ts h\Yh\to\Ug\!$ \text{can be defined by}
$$
H^{(s)}_{ij}+h\Yh\,\mapsto\,-\,F_{ji}^{(s-1)}\cdot(-1)^\bj\ ,
\Tag{3.55555666666}
$$
see \(3.56). Let us extend the comultuplication $\De$ to $\Yh$ by
$\CC[[h]]$-linearity.
The definition \(3.7) implies that
the assignment \(3.55555666666)
defines an isomorphism of Hopf superalgebras.
Let $\psi:\Yh\to\Ug$ be the composition of the projection
$\Yh\to\Yh/\ts h\Yh$ with the isomorphism \(3.55555666666). 

Now let us consider the co-supercommutator
$\ph:\g\to\g^{\ot2}\subset\glN[u]\ot\glN[v]$ determined by \(1.3), where
according to \(1.4) we put
$$
r(u,v)\,=\,
\sum_{ij}\ts E_{ij}\ot E_{ji}\cdot\frac{{(-1)}^{\ts\bj}}{u-v}\ts
+\sum_{ij}\ts E_{ij}\ot E_{-j,\ts -i}\cdot\frac{{(-1)}^{\ts\bj}}{u+v}\ .
\Tag{2.22}
$$
Extend $\ph$ to the co\ts-Poisson bracket $\Ug\to\Ug^{\ot2}$.
Denote this extension by the same letter $\ph$.
Further, denote by $\Den$ the composition of the comultiplication $\De$
on $\Yh$ with the involutive automorphism $\theta$ of the
algebra $\Yh^{\ot2}$, defined in the beginning of Section 1.
To show that $\Yh$ is a quantisation of the
\text{co\ts-Poisson} Hopf superalgebra
$\Ug$ it remains to prove the following proposition.

\proclaim{Proposition 2.4}
For any element $X\in\Yh$ we have the equality
$$
(\psi\ot\psi)\bigl(\ts(\De(X)-\Den(X))/h\ts\bigr)=
\ph\bigl(\ts\psi(X)\bigr)\,.
\Tag{3.555}
$$
\endproclaim

\demo{Proof}
If suffices to verify the equality \(3.555) for the generators $H^{(s)}_{ij}$
of the algebra $\Yh$. By the definitions \(3.55555666666) and \(1.3),\(2.22)
for $s\ge1$ we have in $\Ug^{\ot2}$
$$
\gather
\ph\bigl(\ts\psi\bigl(H^{(s)}_{ij}\bigr)\bigr)=
-\,\ph\bigl(F^{(s-1)}_{ji}\bigr)\cdot(-1)^\bj=
\\
\sum_{1\le r\le s-1}
\left(\,
F_{ki}^{(r-1)}\ot F_{jk}^{(s-r-1)}
\cdot(-1)^{\ts(\bi+\bk+1)(\bj+\bk)}
\,-\,
F_{jk}^{(r-1)}\ot F_{ki}^{(s-r-1)}
\cdot(-1)^{\ts\bj+\bk}
\,\right).
\endgather
$$
On the other hand, by the definition \(3.7) for any $s\ge1$ we have
in $\Yh^{\ot2}$
$$
\align
\De\bigl(H^{(s)}_{ij}\bigr)
&=
H^{(s)}_{ij}\ot1+1\ot H^{(s)}_{ij}+
\sum_{1\le r\le s-1}
h\cdot H^{(r)}_{ik}\ot H^{(s-r)}_{kj}\cdot(-1)^{\ts(\bi+\bk)(\bj+\bk)}\,,
\\
\De^{\hskip-2pt\o}\bigl(H^{(s)}_{ij}\bigr)
&=
H^{(s)}_{ij}\ot1+1\ot H^{(s)}_{ij}+
\sum_{1\le r\le s-1}
h\cdot H^{(r)}_{kj}\ot H^{(s-r)}_{ik}\,.
\endalign
$$
Thus using again the definition \(3.55555666666) 
we get the equality \(3.555) for $X=H^{(s)}_{ij}$
\enddemos 


\section{3. Centre of the Yangian}

\noindent
In this section we will give a description of the centre of the
$\ZZ_2$-graded algebra $\YN$. By definition an element of $\YN$ is
central if it supercommutes with any element of $\YN$.
However, we will see that the centre of $\YN$ consists of even elements only.
We will use some arguments from \hbox{[\ts MNO,\ts Proposition 2.12\ts].}

\newpage 

Let $\tau$ be the antiautomorphism of the $\ZZ_2$-graded algebra $\EndCN$
defined by the assignment
$$E_{ij}\mapsto E_{ji}\cdot{(-1)}^{\bi\ts(\ts\bj\ts+\ts1)}$$
for any $i$ and $j$. Introduce the element of
the algebra $\EndCN^{\ts\ot\ts2}$
$$
Q=\id\ot\tau\ts(\ts P\ts)=
\sum_{i,j}\ts E_{ij}\ot E_{ij}\cdot{(-1)}^{\bi\ts\bj}\ts.
$$
Denote
$$
\bar T(u)=\tau\ts\ot\id\ts\bigl(\ts\widetilde T(u)\bigr)
\in\EndCN\ot\YN[[\ts u\1\ts]].
$$
The following construction of central elements in $\YN$ goes back to
\hbox{[N1,\ts Section \!1].}

\proclaim{Proposition 3.1}
For a certain element $Z(u)\in\YN[[\ts u\1\ts]]$ we have the equality
$$
\bigl(Q\ot1\bigr)\cdot T_1(u)\,\bar T_2(u)=Q\ot Z(u)
\Tag{4.1}
$$
in the algebra $\EndCN^{\ot\ts2}\ot\YN[[\ts u\1\ts]]$. The coefficients
of the series $Z(u)$ are of $\ZZ_2$-degree zero and
belong to the centre of the algebra $\YN$.
\endproclaim

\demo{Proof}
Introduce the rational function
$\bar R\ts(u,v)=\id\ot\tau\ts\bigl(R\ts(u,v)\bigr)$
valued in the algebra $\EndCN^{\ot\ts2}$. One can directly verify
the identity
$$
\bar R\ts(u,v)\,\bar R\ts(-\ts u,-\ts v)=1\ts.
$$
By making use of this identity we derive from \(3.3) the relation
$$
\bigl(\bar R\ts(-\ts u,-\ts v)\ot1\bigr)\cdot\ts T_1(u)\,\bar T_2(v)=
\bar T_2(v)\,T_1(u)\cdot\bigl(\bar R\ts(-\ts u,- \ts v)\ot1\bigr)\ts.
\Tag{4.2}
$$
Let us multiply each side of this relation by $u-v$ and then put $u=v$.
We~obtain 
$$
\bigl(Q\ot1\bigr)\cdot\ts T_1(u)\,\bar T_2(u)=
\bar T_2(u)\,T_1(u)\cdot\bigl(Q\ot1\bigr)\ts.
$$
Since the image of the endomorphism $Q\in\EndCN^{\ts\ot\ts2}$
has dimension one, we get the first statement of Proposition 3.1.
Since $Q$ has $\ZZ_2$-degree zero, the equality \(4.1) shows that
every coefficient of the series $Z(u)$ has $\ZZ_2$-degree zero in $\YN$.

Let us now work with the algebra $\EndCN^{\ot\ts3}\ot\YN[[\ts u\1,v\1\ts]]$.
Using the relations \(3.3),\(4.2) and the definition \(4.1)
we get the equalities
$$
\align
\bigl(\ts Q_{23}\,\bar R_{13}(-\ts u,-\ts v)\,R_{12}(u,v)\ot1\bigr)\cdot
T_1(u)\,T_2(v)\,\bar T_3(v)
&=
\Tag{4.3}
\\
\bigl(\ts Q_{23}\ot1\bigr)\cdot T_2(v)\,\bar T_3(v)\,T_1(u)\cdot
\bigl(\bar R_{13}(-\ts u,-\ts v)\,R_{12}(u,v)\ot1\bigr)
&=
\\
\bigl(\ts Q_{23}\ot Z(v)\bigr)\cdot T_1(u)\cdot
\bigl(\bar R_{13}(-\ts u,-\ts v)\,R_{12}(u,v)\ot1\bigr)
&\,.
\endalign
$$
On the other hand, by \(3.52) we have the identity in 
$\EndCN^{\ot\ts3}(\ts u,v\ts)$
$$
R_{13}(-\ts u,-\ts v)\,P_{23}\,R_{12}(u,v)=P_{23}\cdot
\left(1-\frac1{(u-v)^2}-\frac1{(u+v)^2}\right).
$$
\vskip-10pt
\line{So\hfill}
\vskip-10pt
$$
Q_{23}\,\bar R_{13}(-\ts u,-\ts v)\,R_{12}(u,v)=Q_{23}\cdot
\left(1-\frac1{(u-v)^2}-\frac1{(u+v)^2}\right).
$$
Due to the latter identity we obtain from \(4.3) the equality
$$
T_1(u)\cdot Q_{23}\ot Z(v)=Q_{23}\ot Z(v)\cdot T_1(u).
$$
Hence every coefficient of the series $Z(v)$ commutes with any generator
$T_{ij}^{(s)}$ of the algebra $\YN$
\enddemos

\newpage 

Let us consider the square $S^{\ts2}$ of the antipodal map.
It is an automorphism of the $\ZZ_2$-graded algebra $\YN$.
Here is an alternative definition of the series $Z(u)$.

\proclaim{Proposition 3.2}
\hbox{\!We have $S^{\ts2}\bigl(\ts T_{ij}(u)\bigr)=
T_{ij}(u)\cdot {Z}^{\ts-1}(u)$ for all indices $i\!$ and $\!j$.}
\endproclaim

\demo{Proof}
Definition \(4.1) is equvalent to the collection of relations in
$\!\YN[[\ts u\1]]\!$ 
$$
\sum_i\ T_{ij}(u)\,\widetilde T_{ki}(u)=Z(u)\,\de_{jk}\ts.
\Tag{4.4}
$$
On the other hand, by the definition of the anipodal map $S$ we
have the relations
$$
\sum_i\ T_{ki}(u)\,S\bigl(T_{ij}(u)\bigr)\cdot
{(-1)}^{\ts(\ts\bi\ts+\ts\bk\ts)(\ts\bi\ts+\ts\bj\ts)}=
\de_{jk}\ts.
\Tag{4.5}
$$
By applying the antiautomorphism $S$ to each side of the latter equality we get
$$
\sum_i\ S^{\ts2}\bigl(T_{ij}(u)\bigr)\,\widetilde T_{ki}(u)=\de_{jk}\ts.
$$
By comparing the last equality with \(4.4) we prove Proposition 3.2
\enddemos 

\proclaim{Corollary 3.3}
We have the equalities of formal series in $u\1$
$$
\De\ts\bigl(Z(u)\bigr)=Z(u)\ot Z(u)\ts,
\quad
\ep\bigl(Z(u)\bigr)=1\ts,
\quad
S\ts\bigl(Z(u)\bigr)={Z}^{\ts-1}(u)\ts.
$$
\endproclaim

\demo{Proof}
Let $\theta$ be the involutive automorphism of the algebra
$\YN\ot\YN$ as defined in the beginning of Section 1.
Since \hbox{$\De\o S=\theta\o(S\ot S)\o\De$} we get
$$
\De\ts\bigl(\widetilde T_{ij}(u)\bigr)=
\sum_{k}\,\,
\widetilde T_{kj}(u)\ot\widetilde T_{ik}(u)
$$
from the definition \(3.7). Now by using \(3.7) again we
obtain the first equality in Corollary 3.3 from \(4.4). 
The second equality follows directly from \(4.4).
To obtain the third equality in Corollary 3.3 apply the the antiautomorphism
$S$ to each side of \(4.4) and then use \(4.5) along with Proposition 3.2
\enddemos

\noindent
Observe that due to the relations \(3.2) we have $Z(-\ts u)=Z(u)$. Thus
$$
Z(u)=1+Z^{(2)}\ts u^{-2}+Z^{(4)}\ts u^{-4}+\ldots
$$
for certain central elements $Z^{(2)},Z^{(4)},\ldots\in\YN$.
We have the following theorem.

\proclaim{Theorem 3.4}
\!Elements $Z^{(2)},Z^{(4)},\ldots$ are free generators
of the \hbox{centre of $\YN$.}
\endproclaim

\noindent
We will present the main steps of the proof as separate propositions.
We will make use of the second ascending filtration on the algebra $\YN$.
Take the element $G^{(s)}_{ij}$ of the $\ZZ$-graded algebra
$\grm\YN$ corresponding to \hbox{generator~$T^{(s)}_{ij}\in\YN$.}
Denote 
$$
G^{(s)}=\sum_i\ts\,G^{(s)}_{ii}\cdot(-1)^\bi\ts.
$$
Note that by the relation \(3.2) here $G^{(s)}=0$ if the number $s$
is even. Theorem~2.3 provides an isomorphism between
$\grm\YN$ and the enveloping algebra $\Ug$ of the Lie superalgebra \(1.0).
In particular, the elements $G^{(1)},G^{(3)},\ldots\in\grm\YN$
are algebraically independent.

\proclaim{Proposition 3.5}
For any index $s=2,4,\ldots$ the element of the algebra $\grm\YN$
corresponding to $Z^{(s)}\in\YN$ is $(s-1)\cdot G^{(s-1)}$. 
\endproclaim

\demo{Proof}
Amongst other relations the collection \(4.2) contains the equality
$$
\bigl[\ts T_{ij}(u),\widetilde T_{ji}(v)\ts\bigr]=
\frac{(-1)^{\bi}}{u-v}\ \sum_k\ \bigl(\ts
T_{kj}(u)\,\widetilde T_{jk}(v)-\widetilde T_{ki}(v)\,T_{ik}(u)
\bigr)
$$
for any indices $i$ and $j$. The square brackets here stand for the
supercommutator. By performing summation in this equality over the
index $i$ we get
$$
\sum_i\ T_{ij}(u)\,\widetilde T_{ji}(v)=
1- \sum_{i,k}\ \widetilde T_{ki}(v)\,T_{ik}(u)
\cdot\frac{(-1)^\bi}{u-v}
$$
By setting $u$ equal to $v$ in the latter equality we obtain due to~\(4.4)
that
$$
Z(v)=1- \sum_{i,k}\ \widetilde T_{ki}(v)\,\dot T_{ik}(v)
\cdot{(-1)^\bi}
\Tag{4.6}
$$
\vskip-10pt
\line{where\hfill}
\vskip-5pt
$$
\dot T_{ik}(v)=-\ts\,T_{ik}^{(1)}\ts v^{-2}-\ts2\,T_{ik}^{(2)}\ts v^{-3}
-\ts\ldots
$$
is the first derivative of the formal series $T_{ik}(v)$ with respect to the
parameter~$v$. By the definition of the second filtration on $\YN$
the element of the $\ZZ$-graded algebra $\grm\YN$ corresponding to the
coefficient at $v^{-s}$  in the expansion of the right hand side of \(4.6) is
$(s-1)\cdot G^{\ts(s-1)}$ whenever $s\ge1$
\enddemos

\noindent
To prove Theorem 3.4 it suffices to show that the elements
$G^{(1)},G^{(3)},\ldots$ generate the centre of 
$\grm\YN$. By Theorem 2.3 this means that for the element
$E\in\qN$ defined by \(2.3.0), the elements $E,E\,u^2,E\,u^4,\ldots\in\g$
generate the centre of the $\ZZ_2$-graded algebra
$\Ug$. 

To prove the latter statement we will consider the
following general situation.
Let $\b$ be an arbitrary finite-dimensional Lie superalgebra.
Let $\om$ be any involutive automorphism of $\b$.
Consider the corresponding twisted polynomial current Lie superalgebra
$$
\h=\bigl\{\,X(u)\in\ts\b[u]\,:\ \om\ts\bigl(X(u)\bigr)=X(-u)\,\bigr\}\ts. 
$$

\proclaim{Proposition 3.6}
Suppose that the centre of the Lie superalgebra $\b$ is trivial.
Then the centre of the universal enveloping algebra $\Uh$ is also trivial.
\endproclaim

\demo{Proof}
We will prove that the adjoint action of $\h$ in the
supercommutative algebra $\Sh$ has only trivial invariant elements.
Choose a homogeneous basis $X_1\lc X_n$ in $\b$ and let
$$
[X_p,X_q\ts]=\sum_r\,c_{pqr}\,X_r
\hskip-10pt
$$
where $c_{pqr}\in\CC$ is a structure constant of $\b$.
We put $\br=0$ if the element $X_r\in\b$ is even
and $\br=1$ if this element is odd. 
Assume that for some $h\le n$ we have
$\om\ts(X_r)=X_r$ when $1\le r\le h$ and $\om\ts(X_r)=-X_r$ when $h<r\le n$.

The elements $X_r\ts t^s$ where $1\le r\le h$ when $s=0,2,\ldots$
and $h<r\le n$ when $s=1,3,\ldots$ form a basis in the Lie superalgebra
$\h$. Let us order the set of the pairs $(s,r)$ here lexicographically:
$$
(0,1)\prec\ldots\prec(0,h)
\prec
(1,h+1)\prec\ldots\prec(1,n)\prec\ldots\,\,.
$$
A basis in the supercommutative algebra $\Sh$ is then formed by
all finite ordered products of the elements $(X_r\ts u^s)^{d}$
over the set of pairs $(s,r)$ 
where $d=0,1,2,\ldots$ when $\br=0$ but $d=0,1$ when $\br=1$.

Let us now fix any $\h$-invariant element $Y\in\Sh$.
Let $m$ be the maximal integer such that $X_r\ts u^{\ts m}$ occurs in $Y$
for some index $r$. Suppose that $m$ is even.
Then the element $Y$ is a finite sum
$$
Y\,=\!\sum_{d_1\ldots\,d_h}Y_{d_1\ldots\,d_h}\cdot
(X_1\ts u^{\ts m})^{d_1}\ldots\,(X_h\ts u^{\ts m})^{d_h}
$$
where any factor $Y_{d_1\ldots\,d_h}\in\Sh$
depends only on elements $X_r\ts u^s\in\h$ with $s<m$. 
This factor is zero if $d_p>1$ for some index $p\le h$ with $\bar p=1$.
By our assumption
$$
\operatorname{ad}(X_q\,u)\cdot Y=0\,;\qquad
q=h+1\lc n\,.
\Tag{4.7}
$$
The minimal component of the left hand side of \(4.7) that depends on
elements $X_r\ts u^{m+1}\in\h$ is the sum over $d_1\lc d_h$
of the products in $\Sh$
$$
Y_{d_1\ldots\,d_h}
\sum_{p\le h}\ 
(X_1 u^{\ts m})^{d_1}\!\ldots
(X_p u^{\ts m})^{d_p-1}\!\ldots
(X_h u^{\ts m})^{d_h}
\!\!\sum_{h<r\le n}\!d_p\ts c_{pqr}(X_r u^{m+1})
\!\cdot\hskip-1pt(-1)^{f_p}
$$
where
$$
f_p=
\sum_{p<s\le h}\ \bar s\ts d_s\ts(\bar q+\br)\ts=
\sum_{p<s\le h}\ \bar s\ts d_s\ts\bar p\ts
$$
if $c_{pqr}\neq0$. That component must be equal to zero.
So for all $q,r=h+1\lc n$
$$
\sum_{d_1\ldots\,d_h}\!
Y_{d_1\ldots\,d_h}
\cdot\sum_{p\le h}\ 
(X_1 u^{\ts m})^{d_1}\!\ldots
(X_p u^{\ts m})^{d_p-1}\!\ldots
(X_h u^{\ts m})^{d_h}\cdot d_p\ts c_{pqr}\ts(-1)^{f_p}=0\,.
$$
Thus for any sequence $d_1\lc d_h\!$ of non-negative
integers such that \text{$d_p\le1\!$ if $\bar p=1$\!}
$$
\sum_{p\le h}\ 
Y_{d_1\ldots\,d_p+1\,\ldots\,d_h}\ts(d_p+1)\,c_{pqr}\cdot(-1)^{f_p}=0\ts;
\ \quad
q,r=h+1\lc n\,.
\Tag{4.75}
$$

By our assumption we have along with \(4.7) the collection of equalities in
$\Sh$
$$
\operatorname{ad}(X_q\,u^2)\cdot Y=0\,;\qquad
q=1\lc h\,.
\Tag{4.8}
$$
By considering the minimal
component of the left hand side of \(4.8) that depends
on elements $X_r\ts u^{m+2}\in\h$ we get along with \(4.75)
the equalities
$$
\sum_{p\le h}\ 
Y_{d_1\ldots\,d_p+1\,\ldots\,d_h}\ts(d_p+1)\,c_{pqr}\cdot(-1)^{f_p}=0\,;
\ \quad
q,r=1\lc h\,.
\Tag{4.85}
$$

Let us now make use of the assumption that the centre of 
the Lie superalgebra $\h$ is trivial.
It implies that the system of linear equations on the variables $z_1\lc z_h$
$$
\sum_{p\le h}\  
\bigl[\ts z_p\ts(d_p+1)\ts X_p\cdot(-1)^{f_p}
,X_q\ts\bigl]\,=\,0\,;
\ \quad
q=1\lc n
$$has only the trivial solution. Rewrite the latter system as
$$
\sum_{p\le h}\  
\ts z_p\ts(d_p+1)\ts c_{pqr}\cdot(-1)^{f_p}=0\ts;
\ \quad
q,r=1\lc n
$$
and compare the result with \(4.75),\(4.85). We obtain that
$Y_{d_1\ldots\,d_p+1\,\ldots\,d_h}=0$ for any index $p=1\lc h$.
The sequence ${d_1\lc d_h}$ here can be chosen arbitrarily.
So we get $Y=Y_{0\ts\lc 0}$. The case when $m$ is odd
can be treated similarly
\enddemos

\noindent
By applying Proposition 3.6 to the quotient Lie superalgebra 
$\b=\glN/\ts\CC\!\cdot\!E$ we complete the proof of Theorem 3.4.
Before closing this section let us consider the images in $\UN$ of the elements
$Z^{(2)},Z^{(4)},\ldots\in\YN$ with respect to the evaluation homomorphism
$\pi_N$. By the definition \(3.5) we obtain from \(4.6) that
$$
\pi_N:\,Z^{(2)}\mapsto\,-\,\sum_{k}\,F_{kk}=-\,2E
$$
\vskip-15pt\nt
and
\vskip-5pt
$$
\pi_N:Z^{(s+2)}\mapsto\,-\!\!\!\!\sum_{k_1\lc k_{s+1}}\!\!\!
F_{\ts k_2 k_1}
\ldots\ts
F_{\ts k_{s+1} k_s}
\ts
F_{\ts k_1 k_{s+1}}
\cdot
(-1)^{\ts\bk_1\ts+\ts\ldots\ts+\ts\bk_s}
$$
for each $s=2,4,\ldots$ where the indices $k,\ts k_1\lc k_{s+1}$
run through $\pm1\lc\pm N$. So the elements
$\pi_N(Z^{(2)}),\pi_N(Z^{(4)}),\ldots$
generate the centre of the algebra $\UN$ 
by [\ts S1\ts,\ts Theorem\,1\ts]. In particular, we
have the following corollary to Theorem 3.4.

\proclaim{Corollary 3.7}
The image of the centre of $\YN$ with respect to the evaluation homomorphism
$\pi_N$ coincides with centre of the $\ZZ_2$-graded algebra $\UN$.
\endproclaim

\noindent
Different construction of a distinguished linear basis
in the centre of the $\ZZ_2$-graded algebra $\UN$ was given
in [\ts N3\ts,\ts Section\,4\ts]. In particular, that construction 
yields a $\qN$-analogue of the classical Capelli identity [C]
for the enveloping algebra $\operatorname{U}(\frak{gl}_N)$.
Results of the next section are 
underlying for that construction, cf.\ [\ts N4\ts,\ts Section\,3\ts]\ts.


\section{4. Double of the Yangian}

\noindent
The general notion of a quantum double was introduced in
[\ts D3\ts,\ts Section 13\ts]. Here we consider the quantum double
of the Yangian $\YN;$ cf. [S] and \text{[\ts BL\ts,\ts Section 3.3\ts]\ts.}
We employ it to define the universal $R$-matrix for
the Hopf superalgebra $\YN.$

Firstly consider a 
complex associative unital $\ZZ_2$-graded algebra $\YNP$ with the countable
set of generators $T^{\ts(-s)}_{ij}$ where $s=1,2,\ts\ldots$ and
$i,j=\pm1,\ts\dots\ts,\pm N$.
The $\ZZ_2$-gradation on the algebra $\YNP$
is determined by setting $\deg T_{ij}^{\ts(-s)}=\bi+\bj$ for each $s\ge1$.
To write down defining relations for these
generators we put
$$
\TB_{ij}(v)=
\de_{ij}\cdot1+T_{ij}^{\ts(-1)}+
T_{ij}^{\ts(-2)}\ts v+T_{ij}^{\ts(-3)}\ts v^2+\ldots\in\YNP[[\ts v\ts]]\,.
\Tag{3.10}
$$
Let us now combine all the series~\(3.10) into the single element
$$
\TB\hskip.5pt(v)=
\sum_{i,j}\,\ts\TB_{ij}(v)\ot E_{ij}\in\YNP\ot\EndCN\ts[[\ts v\ts]]\,.
$$

Further, for any positive integer $n$ 
and each $s=1\lc n$ we will denote
$$
\TB_s(v)=\id\ot\iota_s\ts\bigl(\ts\TB(v)\bigr)
\in\YNP\ot\EndCN^{\ts\ot\ts n}\,[[\ts v\ts]]\,.
\Tag{3.122}
$$
Then the defining relations in $\YNP$ can be written as
$$
\align
\TB_1(u)\ts\TB_2(v)\cdot\bigl(1\ot R\ts(u,v)\bigr)
&=
\bigl(1\ot R\ts(u,v)\bigr)\cdot\ts\TB_2(v)\ts\TB_1(u)\,,
\Tag{3.13}
\\
\id\ot\et\,\bigl(\,\TB(v)\bigr)&=\TB(-v)\,.
\Tag{3.14}
\endalign
$$
After multiplying each side of \(3.13) by $u^2-v^2$ it becomes
a relation in the algebra
$$
\YNP\ot\EndCN^{\ot\ts2}\,[[\ts u,v\ts]]\,.
$$
It is equivalent to the collection of relations in the algebra
$\YNP[[\ts u,v\ts]]$
$$
\align
(u^2-v^2)
\cdot
\bigl[\,\ts\TB_{ij}(u&)\ts,\TB_{kl}(v)\ts\bigr]\cdot
(-1)^{\ts\bi\ts\bj\ts+\ts\bi\ts\bl\ts+\ts\bj\ts\bl}\ =
\Tag{3.11}
\\
(u+v)
\cdot
\bigl(\,\TB_{il}(u)\,\TB_{kj}(v)&-\TB_{il}(v)\,\TB_{kj}(u)\bigr)\,\,\ts\ +
\\
\quad\qquad
(u-v)
\cdot
\bigl(\,\TB_{i,-l}(u)\,\TB_{k,-j}(v)
&-\TB_{-i,l}(v)\,\TB_{-k,j}(u)\bigr)
\cdot
(-1)^{\ts\bi\ts+\ts\bj}
\endalign
$$
for all possible indices $i,j$ and $k,l$. Then \(3.14)
is equivalent to the collection of
$$
\TB_{ij}(-\ts v)=\TB_{-i,\ts -j}(v)\,.
\Tag{3.12}
$$

There is a natural structure of $\ZZ_2$-graded
bialgebra on $\!\YNP$. Due to \(3.13) and \(3.14) we can define a
comultiplication $\De:\YNP\to\YNP\ts\ot\YNP$ by
$$
\TB_{ij}(v)\,\mapsto\,\sum_{k}\,\,
\TB_{ik}(v)\ot\TB_{kj}(v)\cdot
(-1)^{\ts(\ts\bi\ts+\ts\bk\ts)(\ts\bj\ts+\ts\bk\ts)},
\Tag{3.17}
$$
similarly to \(3.7). But here the tensor product is taken over the
subalgebra $\CC[[\ts v\ts]]$.
The counit $\ep:\YNP\to\CC$ is determined
so that
$\ep:\ts T^{(-s)}_{ij}\mapsto0$ for $s\ge1$.
Note that $\YNP$ is a bi-superalgebra but not a Hopf superalgebra.
The antipode
is defined for a completion $\YNS$ of $\YNP$ such that
$\TB(0)\in\YNS\ot\EndCN$ is invertible.
We will construct such a completion later in this section.

There is a canonical bilinear pairing
$\langle\,\ts,\ts\rangle:\YN\times\YNP\to\CC$.
We shall describe the corresponding linear map $\beta:\YN\otimes\YNP\to\CC$.
The latter map will be defined following [\ts RTF,\, Section 2\ts] so that
for all numbers $m,n=0,1,2,\ldots$
$$
\align
\hskip17pt\EndCN^{\ot m}\ot\YN&\ot\YNP\ot\EndCN^{\ot n}
\to
\EndCN^{\ot(m+n)}:
\\
T_1(u_1)\ldots T_m(u_m)&\ot\TB_{1}(v_1)\ldots\TB_n(v_n)
\mapsto\!\!
\prod_{1\le k\le m}^\rightarrow
\Bigl(\prod_{1\le l\le n}^\rightarrow\!R_{k,m+l}(u_k,v_l)\Bigr)
\Tag{3.7777777}
\endalign
$$
under the map $\id\ot\beta\ot\id$.
Here $u_1\lc u_m\ts,v_1\lc v_n$ are independent
variables and the product of the rational functions
$R_{k,m+l}(u_k\ts,v_l)$ should be expaned as a formal power series in 
$u_1^{-1}\lc\ts u_m^{-1},v_{1}\ts\lc\ts v_n\ts$.
In particular, when $m=n=0$ we get
the equality $\langle\ts1,\!1\ts\rangle=1$.
Due to the relations \(3.3),\(3.4) and \(3.13),\(3.14) the consistency
of this definition follows from \(3.45),\(3.46) and \(3.6). The following
lemma describes a basic property of the pairing $\langle\,\ts,\ts\rangle$.

\proclaim{Lemma 4.1}
\!Let $s_1\lc s_m\!$ and $r_1\lc r_n\!$ be any numbers from
$\{1,2,\ldots\}\!$.~Then
$$
\bigl\langle\,
T_{i_1j_1}^{\ts(s_1)}\!\ldots T_{i_mj_m}^{\ts(s_m)}\,,\ts
T_{i_{m+1}j_{m+1}}^{\ts(-r_1)}\!\ldots T_{i_{m+n}j_{m+n}}^{\ts(-r_n)}
\bigr\rangle
\neq0
\ts\ \Rightarrow\ 
s_1+\ldots+s_m\ge r_1+\ldots+r_n
$$
for all $m,n=0,1,2,\ldots$ and any choice of the indices
$i_1,j_1\lc i_{m+n}\ts,j_{m+n}\ts$.
\endproclaim

\demo{Proof}
First suppose that $r_1\lc r_n\ge2$. Then by our definition the value
of the pairing in Lemma 4.1 is up to the factor $\pm1$ the coefficient at
$$
E_{i_1j_1}\!\ot\ldots\ot E_{i_{m+n}j_{m+n}}\cdot\ts
u_1^{-s_1}\ldots u_m^{-s_m}\cdot
v_1^{\ts r_1-1}\ldots v_n^{\ts r_n-1}
\Tag{3.91}
$$
in the expansion of the product in
$\EndCN^{\ot\ts(m+n)}\ts
[[\, u_1^{-1}\lc\ts u_m^{-1},v_1\ts\lc\ts v_n\,]]$
$$
\prod_{1\le k\le m}^\rightarrow
\Bigl(\prod_{1\le l\le n}^\rightarrow\!
\bigl(\ts1-\sum_{s\ge1}\,
\frac{v_l^{s-1}}{u_k^s}\ts
P_{\ts k,m+l}\bigl(1+(-1)^sJ_kJ_{m+l})\ts
\bigr)\Bigr)
$$
where we have used \(3.333). If here the coefficient at \(3.91) 
is non-zero then evidently
$$
s_1+\ldots+s_m\ge r_1+\ldots+r_n\ts.
$$

Now suppose that some of the numbers $r_1\lc r_n$ are equal to 1.
Without loss of generality we will assume that $r_1\lc r_p\ge2$ and
$r_{p+1}\lc r_n=1$ for some $p<n$.
Rewrite the product over the indices $k,l$ 
at the right-hand side of \(3.7777777) as
$$
\prod_{1\le l\le p}^\rightarrow
\Bigl(\,
\prod_{1\le k\le m}^\rightarrow\!R_{k,m+l}(u_k\ts,v_l)
\Bigr)\,\cdot\!\!
\prod_{p<l\le n}^\rightarrow
\Bigl(\,
\prod_{1\le k\le m}^\rightarrow\!R_{k,m+l}(u_k\ts,v_l)
\Bigr)\,.
$$
Now the value
of the pairing in Lemma 4.1 is up to the factor $\pm1$ the coefficient at
\(3.91) in the expansion of the product
$$
\align
\prod_{1\le l\le p}^\rightarrow
&\Bigl(\,
\prod_{1\le k\le m}^\rightarrow\!
\bigl(\ts1-\sum_{s\ge1}\,
\frac{v_l^{s-1}}{u_k^s}\ts
P_{\ts k,m+l}(\ts1+(-1)^sJ_kJ_{m+l}\ts)\ts
\bigr)
\Bigr)\ \times
\Tag{3.925}
\\
\quad
\prod_{p< l\le n}^\rightarrow
&\Bigl(\,
\prod_{1\le k\le m}^\rightarrow\!
\bigl(\ts1-\sum_{s\ge1}\,
\frac{v_l^{s-1}}{u_k^s}\ts
P_{\ts k,m+l}(\ts1+(-1)^sJ_kJ_{m+l}\ts)\ts
\bigr)
-1\Bigr)
\endalign
$$
If here that coefficient is non-zero then
$
s_1+\ldots+s_m\ge r_1+\ldots+r_p+n-p
$
\enddemos

We will equip the algebra $\YNP$ with the
descending $\ZZ$-filtration defined by assigning to the generator
$T_{ij}^{\ts(-s)}$ the degree $s$ for any $s\ge1$.
The corresponding $\ZZ$-graded algebra will be denoted by $\gr\hskip.75pt\YNP$.
The formal completion of the algebra $\YNP$ with respect to this
filtration will be denoted by $\YNS$.
We will extend the comultiplication $\De$ on $\YNP$ to the algebra
$\YNS$, and still denote this extension by $\De$. The image
$\De\bigl(\ts\YNS\bigr)$ lies in the formal completion of
the algebra $\YNP\ot\YNP$ with respect to the descending $\ZZ$-filtration,
defined by assigning to the element $T_{ij}^{\ts(-r)}\ot T_{kl}^{\ts(-s)}$
the degree $r+s$. Indeed, with respect to this filtration
$\De(\ts T_{ij}^{\ts(-r)}\bigr)$ is a finite sum of elements
of degree not less than $r$.

Let $G_{ij}^{\ts(-s)}\!\in\!\YNP$ be the element 
corresponding to the generator $T_{ij}^{\ts(-s)}$ of the algebra $\YNP$.
The algebra $\gr\hskip.75pt\YNP$ inherits $\ZZ_2$-gradation~from~$\YNP$
such that for any $s\ge1$ we have $\deg\ts t_{ij}^{\ts(-s)}=\bi+\bj$.
By the relations \(3.12)
we have $$
G_{-i,-j}^{\ts(-s)}={(-1)}^{s+1}\,G_{ij}^{\ts(-s)},
\ \quad
s\ge1\ts.
\Tag{2.2.3}
$$
Furthermore, we can define a bilinear pairing
$$
\langle\,\ts,\ts\rangle:\,\,\gr\hskip.75pt\YN\times\gr\hskip.75pt\YNP\to\CC
\Tag{3.93}
$$
by making 
$$
\bigl\langle\,
t_{i_1j_1}^{\ts(s_1)}\!\ldots t_{i_mj_m}^{\ts(s_m)}\,,\ts
G_{i_{m+1}j_{m+1}}^{\ts(-r_1)}\!\ldots G_{i_{m+n}j_{m+n}}^{\ts(-r_n)}
\,\bigr\rangle
$$
equal to
$$
\bigl\langle\,
T_{i_1j_1}^{\ts(s_1)}\!\ldots T_{i_mj_m}^{\ts(s_m)}\,,\,
T_{i_{m+1}j_{m+1}}^{\ts(-r_1)}\!\ldots T_{i_{m+n}j_{m+n}}^{\ts(-r_n)}
\,\bigr\rangle
$$
\medskip\nt
if
$
s_1+\ldots+s_m=r_1+\ldots+r_n
$
and equal to zero
otherwise.
Here $m,n\ge0$ and $s_1\lc s_m\ts,r_1\lc r_n\ge1$
while the indices $i_1,j_1\lc i_{m+n}\ts,j_{m+n}$ are arbitrary.
This definition is correct due to Lemma 4.1.
Now for each $s=0,1,2,\ldots$ denote by $\YNs$ and $\YNPs$ the
subspaces of degree $s$ in the $\ZZ$-graded
algebras $\gr\hskip.75pt\YN$ and $\gr\hskip.75pt\YNP$
respectively. 

\proclaim{Lemma 4.2}
Restriction of the pairing \(3.93) to $\YNs\times\YNPs\!$ is not
degenerate for any $s\ge0$.
\endproclaim

\demo{Proof}
Fix any integers $s_1\lc s_m\ts,r_1\lc r_n\ge1$ such that
$$
s_1+\ldots+s_m=s_{m+1}+\ldots+s_n\ts.
$$
Without loss of generality we will assume that
$s_1\ge\ldots\ge s_m$ and  $s_{m+1}\ge\ldots\ge s_n\ts.$
Suppose that $r_1\lc r_p\ge2$ while $r_{p+1}\lc r_n=1$ for some $p\ge 0$.
Now we do not exclude the case $p=n$. Let us consider the coefficient~at
$$
u_1^{-s_1}\ldots u_m^{-s_m}\cdot\,v_1^{\ts r_1-1}\ldots v_n^{\ts r_n-1}
\Tag{3.94}
$$
in the expansion of the product \(3.925) as a series in
$u_1^{-1}\lc\ts u_m^{-1},v_1\ts\lc\ts v_n$.
By our assumptions this coefficient can be non-zero only if $m=n$ and
$s_k=r_k$ for all indices $k=1\lc m$. Suppose that this is the case. 
For $r=1,2,\ldots$ denote by $\S_r$ the segment of the sequence $1\lc m$
consisting of all $k$ such that $s_k=r$.
Then the coefficient at \(3.94) in the expansion of \(3.925) equals 
$$
(-1)^m\cdot\,\prod_{r\ge1}
\,\,\Bigl(\,\,
\sum_g
\prod_{\,l\ts\in\S_r}
P_{\ts g(l),m+l}\bigl(\ts1+(-1)^rJ_{g(l)}J_{m+l}\ts\bigr)\Bigr)\ts
\Tag{3.944}
$$
where the index
$g$ runs through the set of all permutations of the sequence $\S_r$.
Note that the factors in each of the above two products commute.

Choose any basis in the space $\YNs$ consisting
of monomials $t_{i_1j_1}^{(s_1)}\ldots t_{i_mj_m}^{\ts(s_m)}$
such that
$$
s_1\ge\ldots\ge s_m\ge1\ts,
\quad
s_1+\ldots+s_m=s
$$
and
$$
i_k\in\{\ts1\lc N\ts\}\ts,
\quad
j_k\in\{\ts\pm1\lc\pm N\ts\}
$$
for $k=1\lc m$ while the number $m\ge0$ can vary.
The above argument using the expression \(3.944) shows
that for any two elements of this basis 
$$
t_{i_1j_1}^{(s_1)}\ldots t_{i_mj_m}^{\ts(s_m)}
\quad
\text{and}
\quad
t_{i_{m+1}j_{m+1}}^{(r_1)}\ldots t_{i_{m+n}j_{m+n}}^{\ts(r_n)}
$$
the value
$$
\bigl\langle\,
t_{i_1j_1}^{\ts(s_1)}\!\ldots t_{i_mj_m}^{\ts(s_m)}\,,\,
G_{j_{m+1}i_{m+1}}^{\ts(-r_1)}\!\ldots G_{j_{m+n}i_{m+n}}^{\ts(-r_n)}
\,\bigr\rangle
$$
\smallskip\nt
is non-zero only if $m=n$ and for each index $k=1\lc m$ we have the equalities
$$
i_{m+k}=i_k\ts,
\ \
j_{m+k}=j_k\ts,
\ \
r_k=s_k\ts.
$$
In the latter case that value up to the factor $\pm1$ is the product
$a!\,b\ts!\,\ldots$ where $a,b,\ldots$ are multiplicities in the sequence
of the triples
$
(i_1,j_1,s_1)\lc(i_m,j_m,s_m).
$
But the products 
$$
G_{j_1i_1}^{\ts(-s_1)}\ldots G_{j_mi_m}^{\ts(-s_m)}\in\ts\gr\ts\YNP
$$
corresponding to elements of our basis in $\gr_s\YN$
span the space $\YNPs$
\enddemos

\noindent
Take the subalgebra $\g^\prime=u\cdot\g$ in the Lie superalgebra
$\glN[u]$, \text{see definition \(1.0).}
Consider the corresponding universal enveloping algebra $\Ugp$.

\proclaim{Corollary 4.3}
The $\ZZ_2$-graded algebras $\gr\hskip.75pt\YNP$ and $\Ugp$ are isomorphic.
\endproclaim

\demo{Proof}
Consider the elements $F_{ij}^{(s)}$ of the universal
enveloping algebra of
$\glN[u]$ with $s\ge0$, defined by \(3.56)\,. Any relation between
these elements follows from \(2.2.1)\ts,\ts\(2.2.2).
On the other hand, the generators $G_{ij}^{\ts(-s)}$ of the algebra
$\gr\ts\YNP$ with $s\ge1$
satisfy \(2.2.3). Due to \(3.11) they also satisfy the relations
$$
\gather
(-1)^{\ts\bi\ts\bj\ts+\ts\bi\ts\bl\ts+\ts\bj\ts\bl}\cdot
\bigl[\ts G_{ij}^{\ts(s)},G_{kl}^{\ts(r)}\ts\bigr]\ts=\,
\de_{kj}\,G_{il}^{\ts(-s-r)}-\ts
\de_{il}\,G_{kj}^{\ts(-s-r)}\,\ts+
\\
\qquad\qquad
\de_{k,-j}\,G_{-i,\ts l}^{\ts(-s-r)}\cdot(-1)^{\bi\ts+\ts\bj\ts+\ts s}
\,-
\de_{-i,l}\,G_{k,-j}^{\ts(-s-r)}
\cdot(-1)^{\bi\ts+\ts\bj\ts+\ts s}
\endgather
$$
for all $r,s\ge1$.
Therefore one can define a homomorphism of the algebra $\Ugp$ onto
$\gr\hskip.75pt\YNP$ by
$$
u\ts F_{ji}^{(s)}\mapsto-\ts\,G_{ij}^{\ts(s+1)}\cdot{(-1)}^{\ts\bi}.
$$
But Lemma 4.2 implies that the kernel of this homomorphism is trivial
\enddemos

\noindent
We formulate the main property of the pairing
$\langle\,\ts,\ts\rangle$ as the next proposition.

\proclaim{Proposition 4.4}
The bilinear map $\langle\,\ts,\ts\rangle:\YN\times\YNP\to\CC$
is non-degenerate bi-superalgebra pairing.
\endproclaim

\demo{Proof}
Lemma 4.1 and Lemma 4.2 show that the pairing $\langle\,\ts,\ts\rangle$
is non-degenerate.
Due to \(3.7) and \(3.17) the definition \(3.7777777) implies
that for any $X,Y\in\YN$ and $X^\prime,Y^\prime\in\YNP$ we have
$$
\bigl\langle\ts X\ts Y,X^{\prime}\bigr\rangle=
\bigl\langle\ts X\!\ot\!Y,\Delta(X^{\prime})\bigr\rangle
\quad\text{and}\quad
\bigl\langle\ts X,X^{\prime}\ts Y^{\prime}\bigr\rangle=
\bigl\langle\ts\Delta(X),X^{\prime}\!\ot\! Y^{\prime}\bigr\rangle
\Tag{3.18}
$$
where we employ the convention
$$
\bigl\langle\ts X\!\ot\!Y,X^{\prime}\!\ot\!Y^{\prime}\bigr\rangle=
\bigl\langle\ts X,Y\bigr\rangle\ts
\bigl\langle\ts X^{\prime},Y^{\prime}\bigr\rangle
\cdot(-1)^{\ts\deg X^\prime\deg Y}
$$
for the homogeneous elements $X$ and $Y^\prime$. Also by definition
we have $\langle\ts1\ts,1\,\rangle=1$.
Moreover, by setting $n=0$ in
the definition \(3.7777777) we get for any $s_1\lc s_m\ge1$
$$
\bigl\langle\ts
T_{i_1j_1}^{\ts(s_1)}\!\ldots T_{i_mj_m}^{\ts(s_m)}\,,\ts1
\ts\bigr\rangle=0\ts,
\ \quad
m\ge1\ts.
$$
Thus $\langle\ts X\ts,1\,\rangle=\ep(X)$ for the counit $\ep$ on $\YN$.
Furthermore, by setting $m=0$
in Lemma 4.1 we obtain for any $r_1\lc r_n\ge1$
$$
\bigl\langle\ts1\ts,
T_{i_{1}j_{1}}^{\ts(-r_1)}\!\ldots T_{i_{n}j_{n}}^{\ts(-r_n)}
\ts\bigr\rangle=0\ts,
\ \quad
n\ge1\ts.
$$
Therefore $\langle\ts1\ts,X^\prime\,\rangle=\ep(X^\prime)$
for the counit $\ep$ on $\YNP$
\enddemos

\noindent
By Lemma 4.1 the pairing $\YN\times\YNP\to\CC$
extends to $\YN\times\YNS$.  
Let us now choose any linear basis in the vector space $\YN$.
An element of this basis will be denoted by $Y_\al$.
There is a system of vectors 
$Y^\al\in\YNS$ dual to this basis. The formal sum of elements from
$\YNS\ot\YN$,
$$
\R\,=\sum_\al\ Y^\al\ot Y_\al
$$
does not depend on the choice of basis in $\!\YN$.
It is called the {\it universal $R\!$-matrix\,} for the Yangian $\YN$.
The {\it double\,} of the Yangian is an associative complex unital algebra
$\DYN$ which contains $\!\YN$ and $\!\YNP$ as subalgebras. 
Moreover, it is generated by these two subalgebras.
We also impose the relations
$$
\R\cdot\De(X)=\Den\!(X)\cdot\R\,,
\ \ \quad
X\in\YNP
\Tag{3.99}
$$
where $\Den$ is composition of the comultiplications $\Delta$
on $\!\YNP$ with the involutive automorphism $\theta$ of the algebra
$\YNP\ot\YNP$.
Either side of the equality \(3.99) makes sense
as a formal sum of elements from $\YNS\ot\DYN$.

The equalities \(3.18) imply that
for the comultiplications $\Delta$ on $\YN$,$\YNS$
$$
\De\ot\id(\ts\R\ts)=\R_{13}\,\R_{\ts23}\,\ts,
\ \quad
\id\ot\De(\ts\R\ts)=\R_{12}\,\R_{\ts13}
\Tag{3.22222222}
$$
where
$$
\R_{12}\,=\sum_\al\ Y^\al\ot Y_\al\ot 1\ts,
\ \ 
\R_{13}\,=\sum_\al\ Y^\al\ot 1\ot Y_\al\ts,
\ \ 
\R_{23}\,=\sum_\al\ 1\ot Y^\al\ot Y_\al\ts.
$$
It follows from \(3.22222222) that $\R\1=\id\ot S\ts(\R)$
for the antipodal map $S$ on $\YN$.

Let us now regard the parameter $z$ in the definition \(3.66)
as a formal parameter. 
Then we get a representation $\YN\to\EndCN\ts[\ts z\ts]$.
We will denote it by $\rho_z.$
Moreover, by  comparing \(3.6),\(3.46) to \(3.13),\(3.14) respectively
we obtain that
$$
\YNP\ot\EndCN\ts[[\ts v\ts]]
\rightarrow
\EndCN^{\ot\ts2}\ts[[\ts z\1,v\ts]]
:\ 
\TB(v)\mapsto R(z,v)
\Tag{3.666}
$$
determines a representation $\!\YNP\to\EndCN\ts[\ts z\1\ts].$
We will denote it~by~$\rho_z^{\ast}$.
More explicitly, for each index $s\ge1$ we have
$$
\rho_z^\ast:\,\ts
T_{ij}^{(-s)}\,\mapsto\,
-\ts\bigl(\ts E_{ji}\,z^{\ts -s}+E_{-j,\ts -i}\ts(-z)^{\ts -s}\ts\bigr)
\cdot{(-1)}^{\ts\bi}\,.
$$
Therefore we can now extend $\rho_z^{\ast}$ to a representation
$\YNS\to\EndCN\ts[[\ts z\1\ts]]$.

\proclaim{Proposition 4.5}
We have $\rho_z^\ast\ot\id\,(\ts\R\ts)=T(z)\!$ and also
$\id\ot\rho_z\,(\ts\R\ts)=\TB(z)$.
\endproclaim

\demo{Proof}
By the definition of our canonical pairing $\YN\ot\YNP\to\CC$, for any $m\ge0$
the element $\TB(z)\in\YNP\ot\EndCN\ts[[\ts z\ts]]$ has the property that
$$
\align
\EndCN^{\ot\ts m}\ot\YN&\ot\YNP\ot\EndCN
\ts\to\,\EndCN^{\ot\ts(m+1)}\,:
\\
T_1(u_1)\ldots T_m(u_m)&\ot\TB(z)
\,\mapsto
R_{1,m+1}(u_1\ts,z)\ldots R_{m,m+1}(u_{m}\ts,z)
\endalign
$$
under the map $\id\ot\beta\ot\id$.
To get the second equality in Proposition 4.5 it suffuces to show
that the element $\id\ot\rho_z\,(\ts\R\ts)\in\YNS\ot\EndCN\ts[\ts z\ts]$
has the same property. Due to the definition of the element $\R$
the latter property amounts to
$$
\gather
\id\ot\rho_z:\,\ts
\EndCN^{\ot\ts m}\ot\YN
\ts\to\,\EndCN^{\ot\ts m+1}\,:
\\
\,T_1(u_1)\ldots T_m(u_m)
\,\mapsto
R_{1,m+1}(u_1\ts,z)\ldots R_{m,m+1}(u_{m}\ts,z)
\endgather
$$
which holds by \(3.66). Proof of the first equality in 
Proposition 4.5 is similar
\enddemos

\proclaim{Corollary 4.6}
Representations $\rho_z$ of $\YN$ and $\rho_z^\ast$ of $\YNP$ determine
a representation of the algebra $\DYN$ in $\EndCN\ts[\ts z,z\1\ts].$
\endproclaim

\demo{Proof}
According to \(3.99) we have to verify for any $X\in\YNP$
the relation
$$
\bigl(\ts\id\ot\rho_u\ts(\R)\bigr)
\cdot
\Bigl(\ts\id\ot\rho_u^\ast\bigl(\Delta\ts(X)\bigr)\!\Bigr)
=
\Bigl(\ts\id\ot\rho_u^\ast\bigl(\Delta^\prime\ts(X)\bigr)\!\Bigr)
\cdot
\bigl(\id\ot\rho_u\ts(\R)\bigr)
$$
in $\YNS\ot\EndCN\ts[\ts u,u\1\ts]$. It suffices to set here $X=\TB_{ij}(v)$.
Due to the definitions \(3.17) and \(3.666) the collection of the resulting
relations for all indices $i,j$ is exactly the defining relation \(3.13)
\enddemos

\noindent
To write down commutation relations in the algebra $\DYN$ we will use
the tensor product $\EndCN\ot\DYN\ot\EndCN$. There is a natural embedding
of the algebra $\EndCN^{\ot2}$ into this tensor product\ts:
$X\ot Y\mapsto X\ot1\ot Y$ for any elements $X,Y\in\EndCN$. Denote by
$\widehat R(u,v)$ the image of \(3.333) with respect to this embedding.
Then we obtain another corollary to Proposition 4.5.

\proclaim{Corollary 4.7}
In
$\EndCN\ot\ts\DYN\ot\EndCN\ts[[\ts u\1,v\ts]]$ we have
$$
\bigl(\ts T(u)\ot1\bigr)\cdot\widehat R(u,v)\cdot\bigl(\ts1\ot\TB(v)\bigr)
=
\bigl(\ts1\ot\TB(v)\bigr)\cdot\widehat R(u,v)\cdot\bigl(\ts T(u)\ot1\bigr)\ts.
\Tag{3.111111}
$$
\endproclaim

\demo{Proof}
Put $X=\TB_{ij}(v)$ in \(3.99). Apply the homomorphism $\rho_u^\ast\ot\id$
to the resulting equality and use the definition \(3.17).
Then we get the equality
$$
\sum_{k}\,
T(u)\!\cdot\!\bigl(\ts\rho_u^\ast\bigl(\TB_{ik}(v)\bigr)\ot\TB_{kj}(v)\bigr)
\!\cdot\!(-1)^{\ts(\ts\bi\ts+\ts\bk\ts)(\ts\bj\ts+\ts\bk\ts)}
=
\sum_{k}\,
\bigl(\ts\rho_u^\ast\bigl(\TB_{kj}(v)\bigr)\ot\TB_{ik}(v)\bigr)
\!\cdot\!T(u)
$$
in $\EndCN\ot\DYN\ts[[\ts u\1,v\ts]]$
by Proposition 4.5.
Due to the definition \(3.666) the collection of the above equalities
for all indices $i\ts,j$ is equivalent to \(3.111111)
\enddemos

\proclaim{Theorem 4.8}
The relation \(3.111111) 
implies the defining relations \(3.99).
\endproclaim

\demo{Proof}
Let $u_1,u_2,\ldots$ be independent formal parameters.
For each $n=1,2,\ldots$ take the tensor product $\nu^\ast$
of the representations $\rho_z^{\ast}:\YNS\to\EndCN\ts[[\ts z\1\ts]]$
with $z=u_1\lc u_n$. Using our descending $\ZZ$-filtration on 
the algebra $\YNP$ and Corollary 4.3\ts,
we can prove that the kernels of all representations $\nu^\ast$
have zero intersection. The proof is similar to the proof of
Proposition 2.2 and is omitted here. Hence it suffices to derive from
the relation \(3.111111) that for any $X\in\YNP$
$$
\nu^\ast\ot\id\,\bigl(\ts\R\cdot\De (X)\bigr)=
\nu^\ast\ot\id\,\bigl(\ts\Den\!(X)\cdot\R\bigr)\,.
\Tag{3.999}
$$
Let us again use
Proposition 4.5 along with the definition \(3.17). The 
collection of all equalities \(3.999) for $X=\TB_{ij}(v)$ with
various indices $i\ts,j$ can be written as 
the single relation in the algebra
$\EndCN^{\ot n}\ot\ts\DYN\ot\EndCN\ts[[\ts u_1\1\lc u_n\1,v\ts]]$
$$
\gather
\bigl(\ts T_1(u_1)\ldots T_n(u_n)\ot1\bigr)
\cdot
\widehat R_{1,n+1}(u_1,v)\ldots \widehat R_{n,n+1}(u_n,v)
\cdot
\bigl(\ts1\ot\TB(v)\bigr)=
\hskip-20pt
\\
\quad
\bigl(\ts1\ot\TB(v)\bigr)
\cdot
\widehat R_{1,n+1}(u_1,v)\ldots \widehat R_{n,n+1}(u_n,v)
\cdot
\bigl(\ts T_1(u_1)\ldots T_n(u_n)\ot1\bigr)
\hskip-20pt
\Tag{3.9999}
\endgather
$$
where $\widehat R_{1,n+1}(u_1,v)\lc\widehat R_{n,n+1}(u_n,v)$
are respectively the images of the elements
$$
R_{1,n+1}(u_1,v)\lc R_{n,n+1}(u_n,v)\in
\EndCN^{\ot(n+1)}\ts[[\ts u_1\1\lc u_n\1,v\ts]]
$$
under the natural embedding of the latter algebra to the former one.
\text{\!But \!using \!\(3.111111)} repeatedly, we obtain \(3.9999)
\enddemos

\nt
Thus we have proved that the relations \(3.111111) together with
the relations \(3.3),\(3.4) and \(3.13),\(3.14) are defining relations for 
the algebra $\DYN$; cf.\ [\ts KT, Section 2\ts]\ts.
 

\section{5. Representations of the Yangian}

\noindent
Here we construct a wide class of irreducible
representations of the algebra $\YN$, by using irreducible
represenations of a certain less complicated algebra $\An$
where $n=1,2,\ldots$ is arbitrary.
The algebra $\An$ was
introduced in [N2] and called the {\it degenerate affine Sergeev algebra\/,}
in honour of the author of [\ts S1,\,S2\ts]\ts.
This is an analogue of the degenerate affine Hecke algebra, which was
employed in [D2] to construct irreducible representations
of the Yangian $\operatorname{Y}(\frak{gl}_N)$
of the general linear Lie algebra $\frak{gl}_N$.
Results presented in this section were reported for the first time
in the summer of 1991 at the Wigner Symposium in Goslar,
Germany. They were also reported in the autumn of 1992 
at the Symposium on Representation Theory in
\line{Yamagata, Japan.
Non-degenerate affine Sergeev algebra is defined in [JN], cf.\ [O2].}

Consider the crossed product $H_n$ of the symmetric group $S_n$ with the
Clifford algebra over the complex field $\CC$ on $n$ anticommuting generators.
These generators are denoted by $c_1\lc c_n$ and are subjected to the relations
$$
c_p^{\ts2}=-1\ts,\qquad
c_p\ts c_q=-\ts c_q\ts c_p\quad\text{if}\quad p\neq q\ts.
$$
The group $S_n$ acts on the Clifford algebra by permutations
of these $n$ generators. Let $w_{pq}\in S_n$ be the transposition of
two numbers $p\neq q$. \text{\!There is a representation}
$H_n\to\EndCN^{\ot n}$ determined by the assignments
$w_{pq}\mapsto P_{pq}$ and $c_p\mapsto J_p\ts$, see definitions
\(1.5) and \(1.55555). The supercommutant of the image of this representation
in $\EndCN^{\ot n}$
coincides by [\ts S2\ts,\,Theorem 3\ts] with the image of the $n$-th tensor
power of defining representation $\UN\to\EndCN$.
By definition, the complex algebra $\An$ is generated
by the algebra $H_n$ and the pairwise commuting elements $x_1\lc x_n$ with the
following relations:
$$
\align
x_p\ts w_{q,q+1}&=w_{q,q+1}\ts x_p\quad\text{if}\quad p\neq q,q+1\ts;
\Tag{5.1}
\\
x_p\ts w_{p,p+1}&=w_{p,p+1}\ts x_{p+1}-1-c_p\ts c_{p+1}\,;
\\
x_p\ts c_q=c_q\ts x_p&\quad\text{if}\quad p\neq q\ts,
\ \quad
x_p\ts c_p=-\,c_p\ts x_p\,.
\endalign
$$
The algebra $\An$ is $\ZZ_2$-graded so that
$\deg c_p=1$ while
$\deg x_p=\deg w_{pq}=0$.

\proclaim{Proposition 5.1}
\!Let $Y\!$
range over a basis in $H_n$ and let each of \text{$s_1\lc s_n\!$ range}
over the non-negative integers. \!Then the products
$Yx_1^{\ts s_1}\!\ldots\,x_n^{\ts s_n}$
\text{\!form a basis in $\An$.\!}
\endproclaim

\demo{Proof}
For $m=0,1,2,\ldots$ one can define a homomorphism
\text{$\ga_m:\An\to\ H_{m+n}\!$ by}
$$
\ga_m:\ \ 
w_{pq}\mapsto w_{m+p,m+q}\,,\quad
c_{p}\mapsto c_{m+p}\,,\quad
x_p\mapsto
\!\!\!\!\sum_{1\le r<m+p}\!\!\!
(\ts1\!+\!\ts c_{m+p}\ts c_r)\,w_{m+p,r}\,.
$$
This can be verified directly by \(5.1). 
Suppose that $m\ge s_1+\ldots+s_n$. Choose for every $p=1\lc n$
a subsequence $\M_p$ in $1\lc m$ of cardinality $s_p$
so that all these subsequences are disjoint.
Write the image of $x_1^{\ts s_1}\!\ldots\,x_n^{\ts s_n}$
under $\ga_m$ as a linear combination of the elements
$c_r\ldots c_{r^\prime}\ts w\in H_{m+n}$ where 
$1\le r<\ldots<r^{\ts\prime}\le m+n$ and $w\in S_{m+n}$.
Consider the terms in this linear combination where $w$ has the
maximal possible length.
Amongst them we find the term
$$
\prod_{1\le p\le n}\Bigl(\,\prod_{r\in\M_p}^\rightarrow w_{m+p,r}\,\Bigr)
$$
which allows us to restore
the exponents $s_1\lc s_n$ and the basis element $Y\in H_n$
from the image
$\ga_m\bigl(\ts Yx_1^{\ts s_1}\!\ldots\,x_n^{\ts s_n}\bigr)$ uniquely.

By using the relations \(5.1) every element of the algebra $\An$
can be expressed as a finite linear combination of the products
$Yx_1^{\ts s_1}\!\ldots\,x_n^{\ts s_n}$. Now take any such a linear
combination and suppose that for all its terms $m\ge s_1+\ldots+s_n$.
Then the above analysis shows that for all the terms,
the images $\ga_m\bigl(\ts Yx_1^{\ts s_1}\!\ldots\,x_n^{\ts s_n}\bigr)$
are linearly independent in $H_{m+n}$
\enddemos

\noindent
Along with pairwise commuting generators $x_1\lc x_n$
we need the non-commuting generators
$$
y_p=x_p-\!\!\sum_{1\le q<p}\!
(\ts1\!+\!\ts c_p\ts c_q)\,w_{pq}\,;\ \quad
p=1\lc n\,.
$$
Observe that the generators $y_1\lc y_n$ belong to the kernel
of the homomorphism $\ga_0:\An\to H_n$ as defined in the proof
of Proposition 5.1. By using this observation,
$$
\align
w\,y_p\ts w^{-1}=y_{w(p)}\,,&\quad w\in S_n\,;
\Tag{5.2}
\\
y_p\ts c_q=c_q\ts y_p\quad\text{if}\quad p\neq q\ts,&
\ \ \quad
y_p\ts c_p=-\,c_p\ts y_p\,.
\endalign
$$
Relations \(5.1) and relations in the first line of \(5.2)
yield the commutation relations
$$
w_{pq}\,[\ts y_p,y_q]=y_p\!-\!\ts y_q+c_p\ts c_q\ts(\ts y_p\!+\!\ts y_q)
\nopagebreak
$$
for the generators $y_p,y_q$ with arbitrary indices $p,q=1\lc n$.

Now take the tensor product of the $\ZZ_2$-graded algebras
$\EndCN^{\ot n}$ and $\An$. Since the elements $x_1\lc x_n\in\An$
pairwise commute, the assignment
$$
\gather
\EndCN\ot\YN\ts[[u^{-1}]]
\ \rightarrow\ 
\EndCN^{\ot(n+1)}\ot\An\ts[[u^{-1}]]\ts:\ \ 
T(u)
\,\mapsto\,
\\
\prod_{1\le p\le n}^\rightarrow
\Bigl(\,
1-P_{\ts 1,p+1}\ot
\frac{1}{u-x_p}\,\ts+\,P_{\ts 1,p+1}J_1J_{p+1}\ot\frac1{u+x_p}
\,\Bigr)
\Tag{5.3}
\endgather
$$
determines a homomorphism $\YN\to\EndCN^{\ot\ts n}\ot\An$, see
\(3.333) and \(3.666666).
As usual, the fractions $1/(u\pm x_p)$ in \(5.3) should be 
expanded as formal power series in $u\1$. 
The next proposition is a key to our construction,
cf.\ \text{[\ts BGHP\ts,\,Section\,2.1\ts]\ts.}

\proclaim{Proposition 5.2}
a) The difference between the product\/ \(5.3) and the sum
$$
1\,-\!\sum_{1\le p\le n}P_{\ts 1,p+1}\ot\frac1{u-y_p}
\,\ts+\!\sum_{1\le p\le n}P_{\ts 1,p+1}J_1J_{p+1}\ot\frac1{u+y_p}
\Tag{5.4}
$$
belongs to the left ideal in the algebra $\EndCN^{\ot(n+1)}\ot\An\ts[[u\1]]$
generated by all the elements
$1-P_{p+1,q+1}\ot w_{pq}$ and $1-J_{p+1}\ts J_{q+1}\ot c_p\ts c_q$ 
with $p\neq q$.

b) The sum\/ \(5.4) commutes with the elements
$P_{p+1,q+1}\ot w_{pq}$ and $J_{p+1}\ot c_p$.
\endproclaim

\demo{Proof}
Part (b) immediately follows from the relations \(5.2).
To prove (a),
we will use induction on $n$. When $n=1$, the equality $x_1=y_1$
provides the induction base. Suppose that $n>1$ and that Proposition 5.2
is true for $n-1$ instead of~$n$. Then the difference
between \(5.3) and \(5.4) equals
$$
\gather
-\ \Bigl(\,
1\,-\!\sum_{1\le p<n}P_{\ts 1,p+1}\ot\frac1{u-y_p}
\,+\!\sum_{1\le p<n}P_{\ts 1,p+1}J_1J_{p+1}\ot\frac1{u+y_p}
\,\Bigr)\,\times
\\
\Bigl(\
P_{\ts 1,n+1}\ot\frac{1}{u-x_n}-P_{\ts 1,n+1}J_1J_{n+1}\ot\frac1{u+x_n}
\,\Bigr)\,+
\\
\Bigl(\
P_{\ts 1,n+1}\ot\frac{1}{u-y_n}-P_{\ts 1,n+1}J_1J_{n+1}\ot\frac1{u+y_n}
\,\Bigr)\ .\  
\endgather
$$
Up to the terms divisible on the right by $1-P_{p+1,n+1}\ot w_{pn}$ or
$1-J_{p+1}J_{n+1}\ot c_p\ts c_n$ with $1\le p<n$, the above sum equals
$$
\gather
P_{1,n+1}
\ot
\biggl(
\frac{1}{u-y_n}-\frac{1}{u-x_n}
+\!\!
\sum_{1\le p<n}
\frac{1}{u-x_n}
\biggl(
\frac1{u-y_p}+
\frac1{u+y_p}\,c_nc_p\!\ts
\biggr)w_{pn}\!\ts
\biggr)\,+
\\
P_{1,n+1}J_1J_{n+1}
\ot
\biggl(
\frac{1}{u+x_n}-\frac{1}{u+y_n}
+\!\!
\sum_{1\le p<n}
\frac{1}{u+x_n}
\biggl(
\frac1{u+y_p}+
\frac1{u-y_p}\,c_nc_p\!\ts
\biggr)w_{pn}
\biggr)
\endgather
$$
where we have used the fact that $x_n$ commutes with $y_1\lc y_{n-1}$.
In these two lines, the second tensor factors
differ by changing $u$ to $-\ts u$. It suffices to show that
in the first line, the second tensor factor equals zero
in $\An\ts[[u^{-1}]]\ts$.
Multiplying this tensor factor on the left by $u-x_n$, on the right
by $u-y_n$, and using the relations \(5.2) in $\An$ we get the sum
$$
(u-x_n)-(u-y_n)+\!\!
\sum_{1\le p<n}(\,1\!+\!\ts c_n\ts c_p\ts)\ts w_{pn}\,.
\nopagebreak
$$
But this sum equals zero by the definition of the element $y_n\in\An$
\enddemos

Take any representation $\xi:\An\to\End(\ts U)$ where the complex
vector space $U$ is $\ZZ_2$-graded, but not necessary
finite-dimensional.
The algebra $\End(\ts U)$
is then $\ZZ_2$-graded. We assume that the homomorphism~$\xi$ preserves
$\ZZ_2$-gradation. 

Take the tensor product of vector spaces $(\CN)^{\ot\ts n}\ot U$.
We identify the tensor product $\EndCN^{\ot n}\ot\End(\ts U)$
with the algebra $\End\bigl((\CN)^{\ot\ts n}\ot U\ts\bigr)$ so that
$$
(A\ot B)\cdot(a\ot b)=(A\ts a)\ot(B\ts b)\cdot(-1)^{\ts\deg a\,\ts\deg  B}.
$$
for any homogeneous $a\in (\CN)^{\ot\ts n},\ts b\in U$ and
$A\in\EndCN^{\ot\ts n}\ts,B\in\End(\ts U)\!$.
There is an action of the hyperoctahedral group
$S_n\ltimes\ZZ_2^n$
in $(\CN)^{\ot\ts n}\ot U$. The transposition $w_{pq}\in S_n$ acts as
the operator $P_{pq}\ot\xi\ts(w_{pq})$ while the generator of
$p$-th direct factor $\ZZ_2$ in $\ZZ_2^n$ acts as the operator
$J_p\ot\xi\ts(c_p)\cdot\sqrt{-1}$. Denote by $\up$ this action.
Let $V$ be the space of co\ts-invariants with respect to $\up$.
This is the quotient space of 
$(\CN)^{\ot\ts n}\ot U$
with respect to the subspace spanned by the images of all the operators
$\up(g)\!-\!1$ with $g\in S_n\ltimes\ZZ_2^n$.
Each of these operators has $\ZZ_2$-degree zero, therefore
the vector space $V$ inherits $\ZZ_2$-gradation from $(\CN)^{\ot\ts n}\ot U$.

Let us expand the element \(5.4) of $\EndCN^{\ot(n+1)}\ot\An\ts[[u^{-1}]]$ 
relative to the basis of standard matrix units in the first
tensor factor $\EndCN$. Then for $s=0,1,2,\ldots$ the coefficient in \(5.4)
\text{at $E_{ij}\ts u^{-s-1}\in\EndCN\ts[[u^{-1}]]$ is}
$$
\sum_{1\le p\le n}
\bigl(\,
\io_p(E_{ji})+(-1)^s\,\io_p(E_{-j,-i})
\ts\bigr)
\ot y_p^s
\cdot{(-1)}^{\ts\bj\ts+1}\,
\in\,
\EndCN^{\ot n}\ot\An
$$
where $\io_p$ denotes embedding of $\EndCN$ to $\EndCN^{\ot n}$
as the $p$-th tensor factor. All these coefficients commute with
the elements $P_{pq}\ot w_{pq}$ and $J_p\ot c_p$ in
the algebra $\EndCN^{\ot n}\ot\An$ by part (b) of Proposition 5.2.
By part (a), one can now define a representation of the
$\ZZ_2$-graded algebra $\YN$ in $V$ by
assigning to the generator $T_{ij}^{(s+1)}$ with $s\ge0$
the action, induced in $V$ by the operator
$$
\sum_{1\le p\le n}
\bigl(\,
\io_p(E_{ji})+
(-1)^{s}\,\io_p(E_{-j,-i})
\ts\bigr)
\ot\xi\ts(\ts y_p^{\ts s})
\cdot{(-1)}^{\ts\bj\ts+1}\,
\Tag{5.45}
$$
in the space $(\CN)^{\ot\ts n}\ot\ts U$.
Correspondence $U\mapsto V$ is the $\YN$-analogue
of the Drinfeld functor [D2]  
for the Yangian $\operatorname{Y}(\frak{gl}_N)$.
\text{Denote our correspondence by $\Cal F_N$.}

For any positive integer $n^\prime$ consider the tensor
product $\An\ot\Anp$ of $\ZZ_2$-graded algebras.
It is isomorphic to the subalgebra in $\Anpn$ generated by 
transpositions $w_{pq}$ where $1\le p<q\le n$ or $n+1\le p<q\le n+n^\prime$,
along with all the elements $c_p$ and $x_p$ where $1\le p\le n+n^\prime$.
Take any $\ZZ_2$-graded representation
\text{$\xi^{\ts\prime}:\ts\Anp\to\End(\ts U^\prime)$.}
Consider the representation of the algebra $\Anpn$ induced from the
representation of $\An\ot\Anp$ in $U\ot U^\prime$.
The algebra $\Anpn$ acts in the vector space
\text{$\Anpn\ot U\ot U^\prime$}
via left multiplication at the first tensor factor. 
We realise the induced representation in the quotient space
of $\Anpn\ot U\ot U^\prime$ by the following relations: 
for homogeneous
$b\in U\ts,b^\prime\in U^\prime\!\ts,
X\in\Anpn\ts,Y\in\An\ts,Y^\prime\in\Anp\!$  
$$
(XZ)\ot b\ot b^\prime=
X\ot\bigl(\ts\xi\ts(Y)\ts b\ts\bigr)\ot
\bigl(\ts\xi^{\ts\prime}(Y^\prime)\ts b^\prime\ts\bigr)\cdot
(-1)^{\ts\deg b\,\,\ts\deg Y^\prime}
\Tag{5.5}
$$
where $Z$ stands for the image of $Y\ot Y^\prime\in\An\ot\Anp$ in the
algebra $\Anpn$. Let us denote by $U\odot U^\prime$ this quotient space.

Consider the representation of the algebra
$\YN$ in the space $V^\prime=\Cal F_N\bigl(\ts U^\prime\ts\bigr)$.
Determine a representation of $\YN$ in
$V\ot V^\prime$ using the comultiplication \(3.7).

\proclaim{Proposition 5.3}
Representation of the algebra $\YN$ corresponding to 
\text{$U\odot\ts U^\prime$,} is equivalent to the representation of $\YN$ in
$V\ot V^\prime$.
\endproclaim

\demo{Proof}
By definition, $V\ot V^\prime$ is a quotient space of
$(\CN)^{\ot n}\ot U\ot(\CN)^{\ot\ts n^\prime}\ot U^\prime$.
Identify the latter tensor product with 
$(\CN)^{\ot\ts(n+n^\prime)}\ot U\ot U^\prime$ via the linear map
$$
a\ot b\ot a^\prime\ot b^\prime
\,\mapsto\,
a\ot a^\prime\ot b\ot b^\prime
\cdot(-1)^{\ts\deg a^\prime\,\ts\deg  b}\,.
$$
Let $W$ be the quotient space of $(\CN)^{\ot\ts(n+n^\prime)}\ot U\ot U^\prime$
corresponding to \text{$V\ot V^\prime$.}
To determine the action of the algebra $\YN$ in $W$,
we can use the \text{representation}
$\YN\to\EndCN^{\ot(\ts n+n^\prime)}\ot\End(U)\ot\End(U^\prime)$,
with respect to which the element $T(u)\in\EndCN\ot\YN\ts[[u\1]]$
is represented by the product 
$$
\align
\prod_{1\le p\le n}^\rightarrow
&\biggl(\ts
1-P_{\ts 1,p+1}\ot\frac1{u-\xi(x_p)\ot1}
\,\ts+\,
P_{\ts 1,p+1}J_1J_{p+1}\ot\frac1{u+\xi(x_p)\ot1}
\,\biggr)\ \times
\\
\prod_{1\le q\le n^\prime\!}^\rightarrow
&\biggl(\ts
1-P_{\ts 1,n+q+1}\ot\frac1{u-1\ot\xi^\prime(x_{q})}
\,\ts+\,
P_{\ts 1,n+q+1}J_1J_{n+q+1}\ot\frac1{u+1\ot\xi^\prime(x_{q})}
\ts\biggr)
\endalign
$$
in the algebra $\EndCN^{\ot\ts(\ts n+n^\prime+1)}\ot\End(U)\ot\End(U^\prime)$.
\text{Here we used \(3.7).}

The space of the representation of $\YN$ corresponding to
$U\odot U^\prime$, is a quotient space of
$(\CN)^{\ot(\ts n+n^\prime)}\ot\Anpn\ot U\ot U^\prime$. The assignment
$$
a\ot a^\prime\ot b\ot b^\prime
\,\mapsto\,
a\ot a^\prime\ot1\ot b\ot b^\prime
\in
(\CN)^{\ot(\ts n+n^\prime)}\ot\Anpn\ot U\ot U^\prime
$$
induces an isomorphism of $W$ to this quotient. This isomorphism commutes
with the action of the algebra $\YN$, since by \(5.5) for any
$1\le p\le n$ and $1\le q\le n^\prime$,
the actions of $x_p$ and $x_{n+q}$
on the class of $1\ot b\ot b^\prime$ in $U\odot U^\prime$
yield the classes of $1\ot\xi\ts(x_p)\ts b\ot b^\prime$ and
$1\ot b\ot\xi^{\ts\prime}\ts(b^\prime)$ respectively
\enddemos

\noindent
To give an example of the correspondence $\Cal F_N:U\mapsto V$,
consider any
{\it principal series representation} of the algebra $\An$. This is the
representation induced from a character $\chi$ of the commutative
subalgebra in $\An$ generated by $x_1\lc x_n$. Note that this
subalgebra is maximal commutative by [\ts N2\ts,\ts Proposition 3.1\ts]. 
Take character $\chi$ such that
$\chi(x_1)=z_1\,\lc\ts\chi(x_n)=z_n$ in $\CC$.
Due to Proposition~5.1, the space $U_{z_1\ldots\ts z_n}$
of the corresponding principal series representation of $\An$
is identified with algebra $H_n$, which acts on itself via
left multiplication. The action of $x_p\in\An$
is then uniquely determined by the assignment
\text{$1\mapsto z_p$ in the space $H_n$.}

\proclaim{Corollary 5.4}
The representation of the algebra $\YN$ corresponding to 
$U_{z_1\ldots\ts z_n}$ is equivalent to the representation\/ \(3.666666).
\endproclaim

\demo{Proof}
The representations of the algebra $\An$ in $U_{z_1\ldots\ts z_n}$
and $U_{z_1}\odot\ldots\odot U_{z_n}$ are equivalent.
Due to Proposition 5.3 it suffices to consider the case $n=1$.
The space of the representation of $\YN$ corresponding to $U_{z}$ is
the quotient of $\CN\ot H_1$ by the relations
$$a\ot c_1\cdot\sqrt{-1}=(Ja)\ot1\cdot(-1)^{\deg a}$$
for homogeneous $a\in\CN$. Assignment $a\ot1\mapsto a$
induces an isomorphism of this quotient to $\CN$. 
By \(5.45), for any $s\ge0$
the generator $T_{ij}^{(s+1)}\in\YN$ acts on the vector $a\ot1\in\CN\ot U_z$
as the operator
$\bigl(\ts E_{ji}\ot z^{\ts s}+E_{-j,\ts -i}\ot (-\ts z)^s
\ts\bigr)\cdot{(-1)}^{\ts\bj\ts+1}\,.$
Comparing this with the definition \(2.19), we complete the proof
\enddemos

\noindent
Let us use
the notion of a $\ZZ_2$-graded irreducibility.
When the $\ZZ_2$-graded vector space $U$ is finite-dimensional,
the representation $\xi$ in $U$ will be called {\it irreducible}
if any $\ZZ_2$-graded subspace in $U$ preserved by $\xi$ is either the
zero space or $U$ itself.

\proclaim{Theorem 5.5}
Suppose that the representation of the $\ZZ_2$-graded algebra $\An$~in~$U$
is finite-dimensional and irreducible. Then the finite-dimensional
representation of the $\ZZ_2$-graded algebra $\YN$ in $V$ is also irreducible. 
\endproclaim

\demo{Proof}
We extend the arguments from \text{[\ts A\ts,\ts Section 4\ts]\ts.}
The algebra $\YN$ contains the 
enveloping algebra $\UN$ as a subalgebra,
see \(3.51). Representation of this subalgebra in
$(\CN)^{\ot n}\ot U$ is the tensor product of 
$n$ copies of 
the defining representation in $\CN$ and the trivial
representation in $U$, see \(5.3).
Let $V_0$ be any non-zero $\ZZ_2$-graded subspace in
$V=\Cal F_N\bigl(\ts U\ts\bigr)$,
preserved by the action of $\YN$. In particular,
$V_0$ is preserved by the action of $\UN$.
By \text{[\ts S2\ts,\,Theorem 3\ts]} there is a \text{$\ZZ_2$-graded}
subspace $U_0\subset U$ preserved by $\xi(H_n)$,
such that $V_0\subset V$ corresponds to $U_0$.
Assume that for any non-zero vector
$b\in U_0$ the image in $V_0$ of the subspace
$(\CN)^{\ot n}\ot b\subset(\CN)^{\ot n}\ot U_0$ is not zero.
Since $\xi$ is irreducible, 
$\xi(\An)\cdot U_0=U\!$.
\text{Let us show that the subspace
$V_0\subset V$ is also $\YN$-\ts cyclic.}

Consider the representation of the algebra $\An$ induced from
the representation of its subalgebra $H_n$ in $U_0$;
$\xi$ is a quotient of this induced representation.
Thus instead of $\xi$ it suffices to take the induced representation.
Realise it in the quotient space $U^\prime$ of $\An\ot U_0$ 
with respect to the relations $(XY)\ot b=X\ot\bigl(\ts\xi\ts(Y)\ts b\ts\bigr)$
for all homogeneous $X\in\An\ts,Y\in H_n$ and $b\in U_0$.
Instead of the subspace $U_0\subset U$ 
it suffices to take
the image 
of the subspace
$1\ot U_0\subset\An\ot U_0$ in this quotient.
Then we have to prove that the subspace 
in $V^\prime=\Cal F_N\bigl(\ts U^\prime\ts\bigr)$ corresponding to
the 
image in $U^\prime$ of $1\ot U_0$, is $\YN$-\ts cyclic.

There is an ascending $\ZZ$-filtration on the algebra $\An$ such that
any generator~$x_p$ is of degree one while $w_{pq}$ and
$c_p$ are of degree zero, see \(5.1). The filtration on $\An$
induces an ascending $\ZZ$-filtration
on the vector space $(\CN)^{\ot n}\ot\An\ot U_0$
and on its quotient $V^\prime$. This filtration on $V^\prime$
is compatible with the action of the algebra $\YN$, when it is
$\ZZ$-filtered so that the degree of the generator $T_{ij}^{\ts(s+1)}$~is~$s$.
But the corresponding $\ZZ$-graded algebra is isomorphic to $\Ug$,
see Theorem 2.3. The corresponding $\ZZ$-graded action of $\Ug$
can be realised in the space of \text{co\ts-invariants} under the action
of the group $S_n\ltimes\ZZ_2^n$ in 
$(\CN)^{\ot n}\ot U_0\ot\CC\ts[\ts x_1\lc x_n]$. 
Here the action in $(\CN)^{\ot n}\ot U_0$
is determined by $\up$ while the action in $\CC\ts[\ts x_1\lc x_n]$
is standard: any permutation $w\in S_n$ acts as $x_p\mapsto x_{w(p)}$,
the generator of the $q$-th factor $\ZZ_2$ in $\ZZ_2^n$ acts as
$x_p\mapsto(-1)^{\de_{pq}\ts}x_p$. Let $W$ be this space of
\text{co\ts-invariants.} The $\ZZ$-graded action of the algebra $\Ug$ in $W$
is induced by its action in the space
$(\CN)^{\ot n}\ot U_0\ot\CC\ts[\ts x_1\lc x_n]$, where the generator
$F_{ij}^{(s)}\in\g$ \text{acts as the operator}
$$
\sum_{1\le p\le n}
\bigl(\ts\io_p(E_{ij})+(-1)^s\,\io_p(E_{-i,-j})\ts\bigr)
\ot1\ot x_p^{\ts s}\ .
$$
We have to prove that the subspace $V_0\ot 1\subset W$
is cyclic under the action of $\Ug$.

Let $u_1\lc u_n$ be complex variables. Let
$\varpi_n$ be the supersymmetrisation map in the space
$$
\bigl(\ts\!\EndCN[u]\ts\bigr)^{\ot n}=\,\EndCN^{\ot n}[u_1\lc u_n]\,,
\hskip-10pt
$$
as in the proof of Proposition 2.2. We have normalised this map
so that $\varpi_n^{\ts2}=\varpi_n$.
Consider the homomorphism $\Ug\to\EndCN^{\ot n}[u_1\lc u_n]$
determined by
$$
F_{ij}^{(s)}\,\mapsto\! 
\sum_{1\le p\le n}
\bigl(\ts\io_p(E_{ij})+(-1)^s\,\io_p(E_{-i,-j})\ts\bigr)\,u_p^{s}\,.
\Tag{5.7}
$$
Image of $\Ug$ under this homomorphism consists of all
polynomials $F(u_1\lc u_n)$ valued in $\EndCN^{\ot n}$ which are
$\varpi_n$-\ts invariant and for each $p=1\lc n$ satisfy
$$
\bigl(\,\id^{\ot\ts (p-1)}\ot\eta\ot\id^{\ot\ts(n-p)}\ts\bigr)\,
F(\ts u_1\lc u_n)=F(\ts u_1\lc\!\!-\!u_p\ts\lc u_n)\,.
$$
This follows from the Poincar\'e\,-Birkhoff\,-Witt theorem 
for Lie superalgebras.

Now consider the subspace in
$(\CN)^{\ot n}\ot U_0\ts\ot\CC\ts[\ts x_1\lc x_n]$ 
consisting of all invariants under the action of the group
$S_n\ltimes\ZZ_2^n$. Denote by $W_{\hskip-1pt\ast}$ this subspace.
Also consider the subspace $W_0$ in the tensor product
$(\CN)^{\ot n}\ot U_0$ consisting of all $\up$-\ts invariants.
We shall prove that the subspace $W_0$ is cyclic 
under the action of the algebra $\EndCN^{\ot n}$ in the first tensor factor.
The above description of the image of $\Ug$ under \(5.7) will then imply, that
the subspace $W_0\ot1\subset W_{\hskip-1pt\ast}$ is \text{$\Ug$-\ts cyclic.}
But this will yield $\Ug$-\ts cyclicity of the subspace
$V_0\ot1\subset W$.

Take any $\up$-invariant inner product $\langle\,\ts,\ts\rangle$
on the vector space $(\CN)^{\ot n}\ot U_0$.
Now suppose that the subspace $W_0\subset(\CN)^{\ot n}\ot U_0$ is not
$\EndCN^{\ot n}$-\ts cyclic. Then we have
$
\bigl\langle\ts(\CN)^{\ot n}\ot b\,,W_0\ts\bigr\rangle=\{0\}
$
for some non-zero vector $b\in U_0$. But this contradicts to our
initial choice of the subspace $U_0\subset U$
\enddemos

\nt
A method for constructing
the irreducible finite-dimensional representations of the 
\text{algebra} $\An$ was developed in [N2].
Restriction of any such representation $U$
to the subalgebra $H_n\subset\An$ is a quotient of the
left regular representation of $H_n$. The restriction of the corresponding
representation $V$ of $\YN$ to the subalgebra $\UN$ is then
a quotient of the representation of $\UN$ in $(\CN)^{\ot n}$. 
But there are irreducible finite-dimensional representations
of the algebra $\UN$, which do not appear as quotients of
the representation in $(\CN)^{\ot n}$ for any $n$\,, see [P]. 
Thus our correspondence $U\mapsto V$ cannot provide all irreducible
representations of the algebra $\YN$.
It would be interesting to give a parametrisation of all irreducible
finite-dimensional representations of the algebra $\YN$; 
\text{ cf.\ [\ts D4\ts,\ts Theorem\,2\ts] and [M]}.


\kern15pt
\line{\bf Acknowledgements\hfill}\section{\,}
\kern-20pt

\noindent
I am grateful to I.\,Cherednik, G.\,Olshanski and A.\,Sudbery
for helpful discussions. I am also grateful to V.\,Drinfeld, P.\,Kulish
and E.\,Sklyanin for their kind interest in this work.
Support from the EPSRC by an Advanced Research Fellowship,
and from the EC under the TMR grant FMRX-CT97-0100, is gratefully acknowledged.


\kern15pt
\line{\bf References\hfill}\section{\,}
\kern-20pt

\itemitem{[A]}
T.\,Arakawa,
{\it Drinfeld functor and finite-dimensional representations of the Yangian},
Commun.\,Math.\,Phys.\ (1999)\ts, 
{\tt math/9807144}.

\itemitem{[A1]}
J.\,Avan,  
{\it Graded Lie algebras in the Yang\ts-Baxter equation},
Phys.\,Lett.\ 
{\bf B\,245}
(1990),
491--496.

\itemitem{[A2]}
J.\,Avan,  
{\it Current algebra realization of $R\!$-matrices associated to $Z_2\!$-graded
Lie algebras},
Phys.\,Lett.\ 
{\bf B\,252}
(1990),
230--236.

\itemitem{[BD]}
A.\,Belavin and V.\,Drinfeld,
{\it The classical Yang\ts-Baxter equation for simple Lie algebras},
Funct.\,Anal.\,Appl.\
{\bf 17}
(1983),
220--221.

\itemitem{[BL]}
{D.\,Bernard and A.\,LeClair},
{\it The quantum double in integrable quantum field theories},
{Nucl.\,Phys.}
{\bf B\,399}
(1993),
709--748. 

\itemitem{[BGHP]}
{D.\,Bernard, M.\,Gaudin, F.\,Haldane and V.\,Pasquier},
{\it Yang-Baxter equation in long-range interacting systems},
J.\,Phys.\
{\bf A\,26}
(1993),
5219--5236.

\itemitem{[C]}
{A.\,Capelli},  
{\it Sur les op\'erations dans la th\'eorie des formes alg\'ebriques},
Math. Ann.\
{\bf 37}
(1890),
1--37.


\itemitem{[D1]}
V.\,Drinfeld,
{\it Hopf algebras and the quantum Yang\ts-Baxter equation},
Soviet Math.\,Dokl.\
{\bf 32}
(1985),
254--258.

\itemitem{[D2]}
{V.\,Drinfeld},
{\it \!Degenerate affine Hecke algebras and Yangians},
Funct.\,Anal.\,Appl.\
{\bf 20}
(1986),
56--58.

\itemitem{[D3]}
{V.\,Drinfeld},
{\it\!Quantum groups},
International Congress of Mathematicians 1986,
Amer.\,Math.\,Soc.\ts,
Providence,
1987,
pp.\ 798--820.

\itemitem{[D4]}
{V.\,Drinfeld},
{\it A \!new realization of Yangians and quantized affine algebras},
\text{Soviet} Math.\,Dokl.\
{\bf 36}
(1988),
212--216.

\itemitem{[FR]}
L.\,Faddeev and N.\,Reshetikhin,
{\it Hamiltonian structures for integrable field theory models},
Theoret.\,Math.\,Phys.\
{\bf 56}
(1983),
847--862.

\itemitem{[JN]}
{A.\,Jones and M.\,Nazarov},
{\it Affine Sergeev algebra and q-analogues of the Young symmetrizers
for projective representations of the symmetric group\,},
J.\,London Math.\,Soc.\
{\bf 78}
(1999),
{\tt q-alg/9712041}. 

\itemitem{[K]}
V.\,Kac,  
{\it Lie superalgebras},
Adv.\,in Math.\ 
{\bf26}
(1977),
8--96.

\itemitem{[KT]}
{S. Khoroshkin and V. Tolstoy},
{\it Yangian double\,},
{Lett.\,Math.\,Phys.}
{\bf 36}
(1996),
373--402.

\itemitem{[LS]}
{D.\,Leites and V.\,Serganova},
{\it Solutions of the classical Yang\ts-Baxter equation
for simple Lie superalgebras},
Theoret.\,Math.\,Phys.\
{\bf 58}
(1984),
16--24.

\itemitem{[M]}
{A.\,Molev},
{\it Finite-dimensional irreducible representations of twisted Yangians},
J.\,Math.\,Phys.\
{\bf 39}
(1998),
5559--5600.

\itemitem{[MM]}
{J.\,Milnor and J.\,Moore},
{\it On the structure of Hopf algebras},
Ann.\,of Math.\
{\bf 81}
(1965),
211--264.

\itemitem{[MNO]}
{A. Molev, M. Nazarov and G. Olshanski},
{\it Yangians and classical Lie algebras},
Russian Math.\,Surveys
{\bf 51}
(1996),
205--282.

\itemitem{[N1]}
M.\,Nazarov,
{\it\!Quantum Berezinian and the classical Capelli identity},
\!Lett.\,Math. Phys.\
{\bf 21}
(1991),
123--131.

\itemitem{[N2]}
M.\,Nazarov,
{\it Young's symmetrizers for projective
representations of the symmetric group},
Adv.\,in Math.\ 
{\bf 127}
(1997),
190--257.

\itemitem{[N3]}
{M.\,Nazarov},
{\it\!Capelli identities for Lie superalgebras},
\!Ann.\ts Scient.\ts\'Ec.\ts Norm.\ts Sup.
{\bf 30}
(1997),
847--872.

\itemitem{[N4]}
{M.\,Nazarov},
{\it Yangians and Capelli identities},
Amer.\,Math.\,Soc.\,Transl.\
{\bf 181}
(1998),
139--163.

\itemitem{[O1]}
{G. Olshanski},
{\it Representations of infinite-dimensional classical groups,
limits of enveloping algebras, and Yangians},
Adv.\,in Soviet Math.\
{\bf 2}
(1991),
1--66.

\itemitem{[O2]}
{G. Olshanski},
{\it Quantized universal enveloping superalgebra
of type $Q$ and a super-extension of the Hecke algebra},
Lett.\,Math.\,Phys.\
{\bf 24}
(1992),
93--102.

\itemitem{[P]}
{I.\,Penkov},
{\it Characters of typical irreducible finite-dimensional
${\frak q}(n)$-modules},
Funct.\,Anal.\,Appl.
{\bf 20}
(1986),
30--37.

\itemitem{[RTF]}
{N.\,Reshetikhin, L.\,Takhtajan and L.\,Faddeev},
{\it Quantization of Lie groups and Lie algebras},
Leningrad Math.\,J.\
{\bf 1}
(1990),
193--225.

\itemitem{[S]}
F.\,Smirnov,
{\it Dynamical symmetries of massive integrable models},
Internat.\,J. Modern.\,Phys.\
{\bf A7}
(1992),
Suppl.\,1B,
813--858.

\itemitem{[S1]}
A.\,Sergeev,
{\it The centre of enveloping algebra for Lie superalgebra $Q(n,\CC)$},
Lett.\,Math.\,Phys.\
{\bf 7}
(1983),
177--179.

\itemitem{[S2]}
A.\,Sergeev,
{\it The tensor algebra of the identity representation
\text{as a module over}
the Lie superalgebras $GL(n,m)$ and $Q(n)$},   
{Math.\,Sbornik}\,\ 
{\bf 51}
(1985),
\text{419--427.\!} 


\bye